\input amstex
\documentstyle{amsppt}
\document
\magnification=1200
\NoBlackBoxes
\nologo
\pageheight{18cm}

\centerline{\bf NEW MODULI SPACES OF POINTED CURVES}

\medskip

\centerline{\bf AND PENCILS OF FLAT CONNECTIONS}

\medskip

\centerline{\bf A.~Losev${}^1$, Yu.~Manin${}^2$}

\medskip

\centerline{\it ${}^1$Institute of Theoretical and Experimental Physics, Moscow, Russia}

\smallskip

\centerline{\it ${}^2$Max--Planck--Institut f\"ur Mathematik, Bonn, Germany}

\bigskip

\hfill{\it Dedicated to William Fulton on the occasion of his 60th birthday}

\bigskip

{\bf Abstract.} It is well known that formal solutions 
to the Associativity Equations are the same as
cyclic algebras over the homology operad 
$(H_*(\overline{M}_{0,n+1}))$ of the moduli
spaces of $n$--pointed stable curves of genus zero. 
In this paper we establish a similar relationship between
the pencils of formal flat connections (or solutions to the
Commutativity Equations) and homology
of a new series $\overline{L}_n$ of pointed stable curves 
of genus zero. Whereas $\overline{M}_{0,n+1}$ parametrizes
trees of $\bold{P}^1$'s with pairwise distinct nonsingular marked points,
$\overline{L}_n$ parametrizes strings of $\bold{P}^1$'s stabilized
by  marked points of two types.
The union of all $\overline{L}_n$'s forms a semigroup
rather than operad, and the role of operadic algebras is taken
over by the representations of the appropriately
twisted homology algebra of this union.

\bigskip

\centerline{\bf 0. Introduction and plan of the paper}

\medskip

One of the remarkable basic results in the theory of the
Associativity Equations (or Frobenius manifolds)
is the fact that their formal solutions are the same as
cyclic algebras over the homology operad 
$(H_*(\overline{M}_{0,n+1}))$ of the moduli
spaces of $n$--pointed stable curves of genus zero. This connection was discovered by physicists, who observed that the data of both types
come from models of topological string theories.
Precise mathematical treatment was given in [KM] and [KMK].

\smallskip

In this paper we establish a similar relationship between
the pencils of formal flat connections (or solutions to the
Commutativity Equations: see 3.1--3.2 below) and homology
of a new series $\overline{L}_n$ of pointed stable curves 
of genus zero. Whereas $\overline{M}_{0,n+1}$ parametrizes
trees of $\bold{P}^1$'s with pairwise distinct nonsingular marked points,
$\overline{L}_n$ parametrizes strings of $\bold{P}^1$'s,
and all marked points with exception of two are allowed
to coincide (see the precise definitions in 1.1 and 2.1).
Moreover, the union of all $\overline{L}_n$'s forms a semigroup
rather than operad, and the role of operadic algebras is taken
over by the representations of the appropriately
twisted homology algebra of this union: see precise definitions in 3.3.

\smallskip

This relationship was discovered on a physical level in [Lo1], [Lo2].
Here we give a mathematical treatment of some of the main
issues raised in these papers.

\smallskip

This paper is structured as follows.

\smallskip

In \S 1 we introduce the notion of $(A,B)$--pointed curves
whose combinatorial structure generalizes that of strings
of projective lines described above. We then describe a construction
of ``adjoining a generic black point'' which allows us
to produce families of such curves and their moduli
stacks inductively. This is a simple variation of one of
the arguments due to F.~Knudsen [1].

\smallskip

In \S 2 we define and study the spaces $\overline{L}_n$
for which we give two complementary constructions. The first one
identifies $\overline{L}_n$ with one of the moduli spaces
of pointed curves. The second one exhibits $\overline{L}_n$
as a well--known toric manifold associated with the
polytope called permutohedron in [Ka2]. These constructions
put $\overline{L}_n$ into two quite different contexts
and suggest generalizations in different directions.

\smallskip

As moduli spaces, $\overline{L}_n$ become components
of the extended modular operad which we define and briefly
discuss in \S 4. We expect that there exists an appropriate
extension of the Gromov--Witten invariants producing
algebras over extended operads involving gravitational
descendants.

\smallskip

As toric varieties, $(\overline{L}_n)$ form one of the several
series related to the generalized flag spaces of classical groups: see [GeSe].
It would be interesting to generalize to other series
our constructions. 

\smallskip

In this paper we use the toric description in order to
prove for $\overline{L}_n$'s an analog of Keel's theorem 
(Theorem 2.7.1) and its extension (Theorem 2.9), crucial
for studying representations of the twisted homology algebra.

\smallskip

This twisted homology algebra $H_*T$ and its relationship
with pencils of formal flat connections are discussed in \S 3,
which contains the main result of this paper: Theorem 3.3.1. 

\bigskip

{\it Acknowledgement.} Yu.~Manin is grateful to M.~Kapranov
who, after having seen the formula $\chi (\overline{L}_n)=n!$,
suggested that $\overline{L}_n$ must be the toric variety
associated with the permutohedron.

\newpage

\centerline{\bf \S 1. $(A,B)$--pointed curves}

\bigskip

\proclaim{\quad 1.1. Definition} Let $A, B$ be two finite disjoint
sets, $S$ a scheme, $g\ge 0$. An $(A,B)$--pointed curve of genus $g$
over $S$ consists of the data
$$
(\pi:\,C\to S;\,x_i:\,S\to C,\ i\in A; \,x_j:\,S\to C,\ j\in B)
\eqno (1.1)
$$
where

\smallskip

(i) $\pi$ is a flat proper morphism whose geometric fibres $C_s$
are reduced and connected curves, with at most ordinary double
points as singularities, and $g=H^1(C_s,\Cal{O}_{C_s}).$

\smallskip

(ii) $x_i, i\in A\cup B,$ are sections of $\pi$ not containing
singular points of geometric fibres.

\smallskip

(iii) $x_i\cap x_j=\emptyset$ if $i\in A,\ j\in A\cup B,\ i\neq j.$

\medskip

Such a curve $(1.1)$ is called stable, if the normalization
of any irreducible component $C^{\prime}$ of a geometric fibre
carries $\ge 3$ pairwise different special points if $C^{\prime}$ is of genus $0$
and $\ge 1$ special points if $C^{\prime}$ is of genus $1.$ 
Special points are inverse images of singular points
and of the structure sections $x_i$.
\endproclaim

\medskip

{\bf 1.2. Remarks.} a) If we put in this definition $B=\emptyset$,
we will get the usual notion of an $A$--pointed (pre)stable curve
whose structure sections are not allowed to intersect
pairwise. Now we divide the sections
into two groups: ``white'' sections $x_i, i\in A$ are not allowed
to intersect any other section, whereas ``black'' sections
$x_j, j\in B$ cannot intersect white ones, but otherwise
are free and can even pairwise coincide. (However,
both types of sections are not allowed to intersect singularities
of fibres).

\smallskip

If we take in this definition one--element set $B=\{*\}$,
we will get a natural bijection between $(A,\{*\})$--pointed curves
and $(A\cup\{*\},\emptyset )$--pointed curves. If $\roman{card}\,B\ge 2,$
the two notions become essentially different.

\smallskip

b) The dual modular graph of a geometric fibre is defined
in the same way as in the usual case (for the conventions
we use see [Ma], III.2). Tails now can be of two
types, and we may refer to them and their marks 
as ``black'' and ``white'' ones as well. Combinatorial
type of a geometric fiber is, by definition,
the isomorphism class of the respective modular
graph with $(A,B)$--marking of its tails.

\smallskip

c) Let $T\to S$ be an arbitrary base change. It produces
from any $(A,B)$--pointed (stable) curve (1.2) over $S$ 
another $(A,B)$--pointed (stable) curve over $T$:
$(C_T;\,x_{i,T}).$

\medskip

{\bf 1.3. A construction.} In this subsection, we start with
an $(A,B)$--pointed curve (1.1) and produce from it another 
$(A,B^{\prime})$--pointed curve: 
$$
(\pi^{\prime} :\,C^{\prime}\to S^{\prime};\,x_i^{\prime},\, i\in A\cup
B^{\prime}).
\eqno (1.2)
$$
The base of the new curve will be $S^{\prime}:=C.$
There will be one extra black mark, say, $*$, so that
$B^{\prime}=B\cup\{*\}$. The new curve and sections
will be produced in two steps. At the first step
we make the base change $C\to S$ as in 1.2 c),
obtaining an $(A,B)$--pointed curve $X:=C\times_S C$,
with sections $x_{i,C}.$
We then add the extra section $\Delta :\,C\to C\times_S C$
which is the relative diagonal, and mark it by $*$.
We did not yet produce an $(A,B^{\prime})$--pointed curve over
$S^{\prime}=C$, because the extra black section can (and generally will)
intersect both singular points of the fibres and
white sections as well. 

\smallskip

At the second step of the construction, we remedy this by birationally
modifying $C\times_S C\to C$ as in [Kn1], Definition 2.3.  
More precisely, we define $C^{\prime}:= \roman{Proj\,Sym}\,\Cal{K}$
as the relative projective spectrum of the symmetric algebra
of the sheaf $\Cal{K}$ on $X=C\times_S C$ defined as
the cokernel of the map
$$
\delta :\,\Cal{O}_X\to  \Cal{J}_{\Delta}\check{}\,\oplus\,
\Cal{O}_X(\sum_{i\in A}x_{i,C}),\ \delta (t) = (t,t).
\eqno(1.3)
$$
Here $\Cal{J}_{\Delta}$ is the $\Cal{O}_X$--ideal of $\Delta$,
and $\Cal{J}_{\Delta}\check{}$ is its dual sheaf
considered as a subsheaf of meromorphic
functions, as in [Kn1], Lemma 2.2 and Appendix.

\smallskip

We claim now that we get an $(A,B^{\prime})$--pointed curve,
because Knudsen's treatment of his modification
can be directly extended to our case. In fact,
the modification we described is nontrivial
only in a neighbourhood of those points, where
$\Delta$ intersects either singular points of the
fibres, or $A$--sections. The $B$--sections do not intersect
these neighborhoods, if they are small enough, and do not
influence the local analysis due to Knudsen ([Kn1], pp. 176--178).

\medskip

{\bf 1.3.1. Remark.} We can try to modify this construction
in order to be able to add an extra white point,
instead of a black one. However, for $\roman{card}\,B\ge 2,$
we will not be able then
to avoid the local analysis of the situation
by referring to [Kn1]. In fact, points where $\Delta$
intersects at least two $B$--sections simultaneously,
will have to be treated anew.

\newpage

\centerline{\bf \S 2. Spaces $\overline{L}_n$}

\medskip

{\bf 2.1. Spaces $\overline{L}_n$.} In this subsection we will inductively
define for any $n\ge 1$ the $(\{0,\infty\}, \{1,\dots ,n\})$--pointed
stable curve of genus zero
$$
(\pi_n: C_n\to \overline{L}_n;\ x_0^{(n)},x_{\infty}^{(n)};\ x_1^{(n)},\dots ,x_n^{(n)}).
\eqno(2.1)
$$
Namely, put 
$$
C_1:=\bold{P}^1,\ \overline{L}_1 = \roman{a\ point},
$$
and choose for $x_0^{(1)},x_{\infty}^{(1)}, x_1^{(1)}$ arbitrary pairwise distinct points.  

\smallskip

If (2.1) is already constructed, we define the next family
$(C_{n+1}\to \overline{L}_{n+1}, \dots )$ as the result
of the application of the construction 1.3 to $C_n/\overline{L}_n.$
In particular, we have a canonical isomorphism $C_n=\overline{L}_{n+1}.$

\smallskip

\proclaim{\quad 2.2. Theorem} a) $\overline{L}_n$ is a smooth separated
irreducible proper manifold of dimension $n-1.$ It represents
the functor which associates with every scheme $T$ the set of the 
isomorphism classes of $(\{0,\infty\}, \{1,\dots ,n\})$--pointed
stable curves of genus zero over $T$ whose geometric fibers
have combinatorial types described below.

\smallskip

The symmetric group $\bold{S}_n$ renumbering the structure sections
acts naturally and compatibly on $\overline{L}_n$ and the universal curve. 
In particular, we can define
the spaces $\overline{L}_B, C_B$ for any finite set $B$, functorial
with respect to the bijections of the sets.

\smallskip

b) Combinatorial types of geometric fibres of $C_n\to\overline{L}_n$
are in a natural bijection with ordered partitions
$$
\{1,\dots ,n\}=\sigma_1\cup\,\dots\,\cup\, \sigma_l,\ 1\le l\le n,\ \sigma_i\ne\emptyset .
\eqno(2.2)
$$
Partition (2.2) corresponds to the linear graph with vertices
$(v_1,\dots ,v_l)$ of genus zero, edges joining $(v_i,v_{i+1}), 1\le i\le l-1$,
$A$--tail $0$ at the vertex $v_1$, $A$--tail $\infty$ at the vertex 
$v_l$, and $B$--tails marked by the elements of $\sigma_i$ at the
vertex $v_i$.

\smallskip

We will call $l=l(\sigma )$ the length of the partition $\sigma$ as in (2.2).

\smallskip

c) Denote by $L_{\sigma}$ the set of all points of $\overline{L}_n$
corresponding to the curves of the combinatorial type $\sigma$,
and by $\overline{L}_{\sigma}$ its Zariski closure.
Then $L_{\sigma}$ are locally closed subsets, and we have
$$
\overline{L}_{\sigma} =\coprod_{\tau \le\sigma} L_{\tau}
\eqno(2.3)
$$
where $\tau \le\sigma$ means that $\tau$ is obtained from
$\sigma$ by replacing each $\sigma_i$ by an ordered partition
of $\sigma_i$ into non--empty subsets.

\smallskip

d) For every $\sigma$, there exists a natural isomorphism
$$
L_{|\sigma_1|}\times \dots \times L_{|\sigma_l|}
\to L_{\sigma}
\eqno(2.4)
$$
such that the pointed curve induced by this isomorphism over
$L_{|\sigma_1|}\times \dots \times L_{|\sigma_l|}$
can be obtained by clutching the curves $C_{|\sigma_i|}/L_{|\sigma_i|}$
in an obvious linear order ($\infty$--section of the $i$--th curve
is identified with the $0$--section of the $(i+1)$--th curve,
see [Kn1], Theorem 3.4), and subsequent remarking of the $B$--sections.

\smallskip

In particular, $L_{\sigma}$ is a smooth 
irreducible submanifold of codimension $l(\sigma )-1.$

\smallskip

The similar statements hold for the closed strata $\overline{L}_{\sigma}.$
\endproclaim
 
\smallskip

{\bf Proof.} Properness and smoothness follow by induction
and Knudsen's local analysis which we already invoked.  

\smallskip

The statement about the combinatorial types is proved
by induction as well. In fact, if everything
is already proved for $C_n$, then we must look at a
geometric fibre $C_{n,s}$ of $C_n$ and see what happens to it
after the blow up described in 1.3. If $\Delta$
intersects a smooth point of $C_{n,s}$, not
coinciding with $x_{0,s}, x_{\infty,s}$, nothing happens,
except that we get a new black point on this fibre,
and a new tail at the respective vertex of the
dual graph. If $\Delta$
intersects an intersection point of two neighboring
components of $C_{n,s}$,
then after blowing up these two components
become disjoint, and we get a new component
intersecting both of them, with a new black point on it.
The linear structure
of the graph is preserved. Finally, if $\Delta$
intersects $C_{n,s}$ at $x_{0,s}$ or $x_{\infty,s}$,
then after blowing up we will get a new end component,
with $x_{0,s}$, resp. $x_{\infty,s}$ and the new
black point on it. Thus the new combinatorial types
will be linear and indexed by partitions of $(n+1)$.
To check that all partitions are obtained in this way,
it suffices to remark that $\Delta$, being the relative diagonal,
can intersect the
fibre of a given type at any point. 

\smallskip

In order to check the statement about the functor represented
by $\overline{L}_n$ we apply the following inductive reasoning. For $n=1$
the statement is almost obvious. In fact, let $\pi :\,C\to S$
be a $(\{0,\infty\}, \{1\})$--pointed
stable curve of genus zero over $T$. From the stability
it follows that all geometric fibres are projective lines.
Since the three structure sections pairwise do not
intersect, the family can be identified with $\bold{P}^1\times T$
endowed with three constant sections. This means that it
is induced by the trivial morphism $T\to\overline{L}_1$.  

\smallskip

Assume that the statement is true for $n$. In order to prove it
for $n+1$, consider a $(\{0,\infty\}, \{1,\dots ,n+1\})$--pointed
stable curve of genus zero $\pi :\,C\to T.$ First of all,
one can produce from it a $(\{0,\infty\}, \{1,\dots ,n\})$--pointed
stable curve of genus zero $\pi :\,C^{\prime}\to T$ obtained
by forgetting $x_{n+1}$ and subsequent stabilization. The
respective map $C^{\prime}\to C$ is given by the relative
projective spectrum of the algebra $\sum_{k=0}^{\infty}
\pi_*(\Cal{K}^{\otimes k})$ where $\Cal{K}:=\omega_{C/T}(x_0+x_1 +
\dots +x_n+x_{\infty})$. By induction, $C^{\prime}$ is induced
by a morphism $p:\,T\to\overline{L}_n$. Addition of an extra black section
to $C^{\prime}$ and subsequent stabilization boils down exactly
to the construction 1.3 applied to $C^{\prime}/T$ which allows us
to lift $p$ to a unique morphism $q:\,T\to\overline{L}_{n+1}.$

\smallskip

Separatedness is checked by the standard deformation arguments.

\smallskip

The statement about renumbering follows from the description
of the functor.

\smallskip

A similar adaptation of Knudsen's arguments allows us
to prove the remaining statements, and we leave them to the reader.

\smallskip

Notice that below we will give another direct description
of the spaces $\overline{L}_B$ and all the structure morphisms connecting them
in terms of toric geometry. This will provide
easy alternate proofs of their properties. Except for
\S 4, we can restrict ourselves to this alternate description.

\medskip

{\bf 2.2.1. Remark.} Dual graphs of the degenerate fibers of $C_n$ over
$\overline{L}_n$ come with a natural orientation from $x_0$ to
$x_{\infty}$. We could have allowed ourselves not to distinguish
between the two white points, interchanging them by isomorphisms, but this would produce several
upleasant consequences. First, our manifolds would become
actual stacks, starting already with $\overline{L}_1$.
Second, we would have lost the toric interpretation of these
spaces. Third, and most important, we would meet an ambiguity
in the definition of the multiplication between the homology spaces:
see (3.5) below. With our choice, we can simply introduce
the involution permuting $x_0$ and $x_{\infty}$ as a part of the structure
and look how it interacts with other parts.

\medskip

\proclaim{\bf \quad 2.3. Theorem} $\overline{L}_n$ has no odd cohomology. 
Let
$$
p_n(q):=\sum_{i=0}^{n-1} \roman{dim}\,H^{2i}(\overline{L}_n)\,q^i
\eqno(2.5)
$$
be the Poincar\'e polynomial of $\overline{L}_n$. Then we have
$$
1+\sum_{n=1}^{\infty} \frac{p_n(q)}{n!}y^n =
\frac{q-1}{q-e^{(q-1)y}} \in \bold{Q}[q][[y]].
\eqno(2.6)
$$
Letting here $q\to 1$ we get $\dfrac{1}{1-y}$ so that
$\chi (\overline{L})=n!.$       
\endproclaim

\smallskip

{\bf Proof.} Since $\overline{L}_n$ are defined over
$\bold{Q}$, we can apply the classical Weil's technique of
counting points over $\bold{F}_q$ (thus treating $q$ not as a formal variable
but as a power of prime). After the counting is done, we will see that
$\roman{card}\,\overline{L}_n(\bold{F}_q)$ is a polynomial in $q$
with positive integer coefficients, so that we can right away 
identify it with $p_n$:
$$
p_n(q)=\roman{card}\,\overline{L}_n(\bold{F}_q)
\eqno(2.7)
$$
The latter number can be calculated by directly applying (2.3)
to the one--element partition $\sigma$, so that we get
$$
\frac{p_n(q)}{n!} =\sum_{l=1}^n \sum\Sb (s_1,\dots ,s_l)\\
s_1+\dots +s_l=n\\s_i\ge 1\endSb
\frac{(q-1)^{s_1-1}}{s_1!}\dots \frac{(q-1)^{s_l-1}}{s_l!} 
$$
$$
\sum_{l=1}^n\left[ \roman{coeff.\ of}\ x^{n-l}\ \roman{in}\
\left( \frac{e^x-1}{x} \right)^l \right]\cdot (q-1)^{n-l}.
$$
Inserting this in the left hand side of (2.6) and summing over $n$
first, we obtain
$$
\sum_{n=1}^{\infty} \frac{p_n(q)}{n!} y^n=
\sum_{l=1}^{\infty}\sum_{n=l}^{\infty}
\left[ \roman{coeff.\ of}\ x^{n}\ \roman{in}\
({e^x-1})^l \right]\cdot (q-1)^{n}
$$
$$
=\sum_{l=1}^{\infty}\frac{1}{(q-1)^l}\,(e^{(q-1)y}-1)^l
$$
which gives (2.6).

\medskip

{\bf 2.3.1. Special cases.} Here is a list of the Poincar\'e
polynomials for small values of $n$:
$$
p_1=1,\ p_2=q+1,\ p_3=q^2+4q+1,\ p_4=q^3+11q^2+11q+1,
$$
$$
p_5=q^4+26q^3+66q^2+26q+1,\ p_6=q^6+57q^5+302q^4+302q^2+57q+1.
$$
The rank of $H^2(L_n)$ is $2^n-n-1$. Individual coefficients of
of $p_n(q)$ are well known in combinatorics. They are called Euler numbers: 
$$
a_{n,i} = \roman{dim}\,H^{2i}(\overline{L}_n) \,.
$$

\medskip

{\bf 2.4. $\overline{L}_n$ and toric actions.} Let $\varepsilon$
be the trivial partition of $B$ of length one. The ``big cell''
$L_{\varepsilon}$ of $\overline{L}_B$ (see 2.2 c)) 
has a canonical structure of
the torsor (principal homogeneous space) over the torus $T_B:=\bold{G}_m^B/\bold{G}_m$
(where the subgroup $\bold{G}_m$ is embedded diagonally).
In fact, $\bold{P}^1\setminus \{x_0,x_{\infty}\}$ is a
$\bold{G}_m$--torsor, and the respective action of
$\bold{G}_m^B$ on $L_{\varepsilon}$, moving $x_i,\,i\in B$
via the $i$--th factor, produces an isomorphic marked curve
exactly via the action of the diagonal. 

\smallskip

Similarly, every stratum $L_{\sigma}$ is a torsor
over $T_{\sigma}:=\prod_i T_{\sigma_i}$ (see (2.4)), and there is a canonical
surjective morphism $T_B\to T_{\sigma}$ so that $L_B$
is a union of $T_B$--orbits. In order to show that
$L_B$ is a toric variety, it remains to show that these
actions are compatible. This again can be done using the
explicit construction of $\overline{L}_n$ and induction.
For a change, we will provide a direct toric construction.
We start with a more systematic treatment of
the combinatorics involved.

\medskip

{\bf 2.4.1. Partitions of finite sets.} For any finite
set $B$, we call a partition $\sigma$ of $B$ {\it a totally
ordered set of non--empty subsets of $B$ whose union is
$B$ and whose pairwise intersectons are empty.}
If a partition consists of $N$ subsets, it is called
$N$--partition. If its components are denoted $\sigma_1,\dots ,\sigma_N$,
or otherwise listed,
this means that they are listed
in their structure order. Another partition can be denoted
$\tau$, $\sigma^{(1)}$ etc. Notice that no particular ordering of $B$
is a part of the structure.
This is why we replaced $\{1,\dots ,n\}$
here by an unstructured set $B.$

\smallskip

Let $\sigma$ be a partition of $B$, $i,j\in B.$
We say that $\sigma$ {\it separates $i$ and $j$} if they
belong to different components of $\sigma$. We then write
$i\sigma j$ in order to indicate that the component
containing $i$ comes earlier that the one containing $j$
in the structure order.

\smallskip

Let $\tau$ be an $N+1$--partition of $B$. If $N\ge 1,$ it determines
a well ordered family of $N$ 2--partitions $\sigma^{(a)}$:  
$$
\sigma^{(a)}_1:=\tau_1\cup\dots\cup\tau_{a},\
\sigma^{(a)}_2:=\tau_{a+1}\cup\dots\cup\tau_{N},\ a=1,\dots ,N\, .
\eqno(2.8)
$$
In reverse direction, call a family of 2--partitions $(\sigma^{(i)})$ 
 {\it good} 
if for any $i\ne j$ we have $\sigma^{(i)}\ne \sigma^{(j)}$
and either  $\sigma^{(i)}_1\subset \sigma^{(j)}_1,$
or $\sigma^{(j)}_1\subset \sigma^{(i)}_1.$ Any good family
is naturally well--ordered by the relation 
$\sigma^{(i)}_1\subset \sigma^{(j)}_1$, and we will
consider this ordering as a part of the structure. If a good family of 2--partitions
consists of $N$ members, we will usually choose
superscripts $1,\dots ,N$ to number these partitions
in such a way that $\sigma^{(i)}_1\subset \sigma^{(j)}_1$
for $i<j.$

\smallskip

Such a good family produces one $(N+1)$--partition $\tau$:
$$
\tau_1:=\sigma_1^{(1)},\ \tau_2:=\sigma_1^{(2)}\setminus
\sigma_1^{(1)},\ \dots ,\  
\tau_N:=\sigma_1^{(N)}\setminus
\sigma_1^{(N-1)},\ \tau_{N+1}=\sigma_2^{(N)}.
\eqno(2.9)
$$ 
This correspondence between good $N$--element families of
2--partitions and $(N+1)$--partitions is one--to--one, because
clearly $\sigma_1^{(i)}=\tau_1\cup\dots\cup\tau_i$ for $1\le i\le N.$

\smallskip

Consider the case when $\tau^{(1)}=\sigma$ is a 2--partition,
and $\tau^{(2)}=\tau$ is an $N$--partition, $N\ge 2$.
Their union is good, iff there exists $a\le N$
and a 2--partition $\alpha =(\tau_{a1},\tau_{a2})$ 
of $\tau_a$ such that
$$
\sigma=(\tau_1\cup\dots\cup\tau_{a-1}\cup
\tau_{a1}, \tau_{a2}\cup
\tau_{a+1}\cup\dots\cup\tau_{N}).
\eqno(2.10)
$$
In this case we denote
$$
\sigma *\tau = \tau (\alpha ):=
(\tau_1,\dots ,\tau_{a-1},
\tau_{a1}, \tau_{a2},
\tau_{a+1},\dots ,\tau_{N}).
\eqno(2.11)
$$

\smallskip

\proclaim{\quad 2.4.2. Lemma}  Let $\tau$ be a
partition of $B$ of length $\ge 1,$ and $\sigma$ a 2--partition.
Then one of the three mutually exclusive cases occurs:

\smallskip

(i) $\sigma$ coincides with one of the partitions
$\sigma^{(a)}$ in (2.8). In this case we will say that
$\sigma$ breaks $\tau$ between $\tau_a$ and $\tau_{a+1}.$

\smallskip

(ii) $\sigma$ coincides with one of the partitions (2.10).
In this case we will say that
$\sigma$ breaks $\tau$ at $\tau_a$.

\smallskip

(iii) None of the above. In this case we will say that
$\sigma$ does not break $\tau$. This happens exactly when there is
a neighboring pair $(\tau_b,\tau_{b+1})$ of elements of $\tau$
with the following property:
$$
\tau_b\setminus\sigma_1\ne \emptyset,\ 
\tau_{b+1}\cap\sigma_1\ne \emptyset .
\eqno(2.12)
$$
We will call $(\tau_b,\tau_{b+1})$ a bad pair (for $\sigma$). 
\endproclaim

\smallskip

{\bf Proof}. Consider the sequence of sets
$$
\sigma_1\cap\tau_1, \sigma_1\cap\tau_2, \dots , \sigma_1\cap\tau_N .
$$
Produce from it a sequence of numbers 0,\,1,\,2 by the following
rule: replace $\sigma_1\cap\tau_b$ by 2, if it coincides
with $\tau_b$, by 0 if it is empty, and by 1 otherwise.
Cases (i) and (ii) above together will furnish all sequences
of the form $(2\dots 20\dots 0)$, $(2\dots 210\dots 0)$,
$(10\dots 0)$. Each remaining admissible sequence
will contain at least one pair of neighbors from the list
01, 02, 11, 12. For the respective pair of sets,
(2.12) will hold.

\medskip

{\bf 2.5. Fan $F_B$.} In this subsection we will describe
a fan $F_B$ in the space $N_B\otimes{\bold{R}}$, where
$N_B:=\roman{Hom}\,(\bold{G}_m,T_B)$, $T_B:=\bold{G}_m^B/\bold{G}_m$
as in the beginning of 2.4. Up to notation, we use [Fu]
as the basic reference on fans and toric varieties.

\smallskip

Clearly, $N_B$ can be canonically identified with
$\bold{Z}^B/\bold{Z}$, the latter subgroup being embedded
diagonally. Similarly, $N_B\otimes{\bold{R}}=\bold{R}^B/\bold{R}$.
We will write the vectors of this space (resp. lattice) as functions
$B\to \bold{R}$ (resp. $B\to \bold{Z}$) considered modulo
constant functions. For a subset $\beta\subset B$, let $\chi_{\beta}$
be the function equal 1 on $\beta$ and 0 elsewhere.

\smallskip

\proclaim{\quad 2.5.1. Definition} The fan $F_B$ consists
of the following $l$--dimensional cones $C(\tau )$ labeled by 
$(l+1)$--partitions $\tau$ of $B$.

\smallskip

If $\tau$ is the trivial 1--partition, $C(\tau )=\{0\}$.

\smallskip

If $\sigma$ is a 2--partition, $C(\sigma )$ is generated
by $\chi_{\sigma_1}$, or, equivalently, $-\chi_{\sigma_2}$,
modulo constants.

\smallskip

Generally, let $\tau$ be an $(l+1)$--partition, and 
$\sigma^{(i)},\,i=1,\dots,l$,
the respective good family of 2--partitions (2.9). Then
$C(\tau )$ as a cone is generated by all $C(\sigma^{(i)})$. 
\endproclaim

\medskip

It is not quite obvious that $F_B$ is well defined.
We sketch the relevant arguments.

\smallskip

First, all cones $C(\tau )$ are strongly convex. In fact,
according to [Fu], p. 14, it suffices to check that
$C(\tau )\cap (-C(\tau ))={0}$. But $C(\tau )$
consists of classes of linear combinations with non--negative coefficients of functions
$$
\chi_{\tau_1},\, \chi_{\tau_1}+\chi_{\tau_2},\,
\dots ,\, \chi_{\tau_1}+\dots +\chi_{\tau_l}
$$
if $\tau$ has length $l+1$. Non--vanishing function
of this type cannot be constant.

\smallskip

Second, the same argument shows that $C(\tau )$ is actually
$l$--dimensional.

\smallskip

Third, since the cone $C(\tau )$ is simplicial, one sees that
$(l-1)$--faces of $C(\tau )$ are exactly $C(\tau^{(i)})$
where $\tau^{(i)}$ is obtained from $\tau$ by uniting
$\tau_i$ with $\tau_{i+1}$, which is equivalent to
omitting $C(\sigma^{(i)})$ from the list of generators.
More generally, $C(\tau^{\prime})$ is a face of
$C(\tau )$ iff $\tau\le\tau^{\prime}$ as in (2.3),
that is, if $\tau$ is a refinement of $\tau^{\prime}$.

\smallskip

Fourth, let $C(\tau^{(i)}),\,i=1,2,$ be two cones.
We have to check that their intersection is a cone
of the same type. An obvious candidate is $C(\tau )$
where $\tau$ is the crudest common refinement of $\tau^{(1)}$
and $\tau^{(2)}$. This is the correct answer. 

\smallskip

In order to see this, let us a give a different description
of $F_B$ which will simultaneously show that
the support of $F_B$ is the whole space. Let $\chi :\bold{B}\to\bold{R}$
represent an element $\bar{\chi}\in N_B\otimes\bold{R}.$ It defines
a unique partition $\tau$ of $B$ consisting of the level sets of $\chi$
ordered in such a way that the values of $\chi$ decrease.
Clearly, $\tau$ depends only on $\bar{\chi}$, and $\chi$ modulo constants
can be expressed as a linear combination of $\chi_{\tau_1}+\dots
+\chi_{\tau_i}$, $1\le i\le l$ with positive coefficients.
In other words, $\chi$ belongs to the interior part of $C(\tau )$.
On the boundary, some of the strict inequalities between 
the consecurive values of $\chi$ become equalities. This proves
the last assertion. 

\smallskip

We see now that $F_B$ satisfies the definition of [Fu], p. 20,
and so is a fan.

\medskip

{\bf 2.6. Toric varieties $\overline{\Cal{L}}_B$.} We now
define $\overline{\Cal{L}}_B$ (later to be identified with
$\overline{L}_B$) as the toric variety
associated with the fan $F_B$.

\smallskip

To check that it is smooth, it suffices to show that
each $C(\tau )$ is generated by a part of a basis of $N_B$
(see [Fu], p. 29). In fact, let us choose a total
ordering of $B$ such that if $i\in\tau_k,\,j\in\tau_l$
and $k<l$, then $i<j$. Let $B_k\subset B$ consist of the first $k$
elements of $B$ in this ordering. Then the classes of the characteristic functions
of $B_1, B_2,\dots ,B_{n-1}$, $n=\roman{card}\,B$, form a basis
of $N_B$, and $\{\chi_{\sigma^{(i)}}\}$ is a part of it.

\smallskip

To check that  $\overline{\Cal{L}}_B$ is proper, we have to show that
the support of $F_B$ is the total space. We have already proved this.

\smallskip

As any toric variety, $\overline{\Cal{L}}_B$  
carries a family of subvarieties which are the closures of the orbits of $T_B$
and which are in a natural bijection with the cones $C(\tau )$ in $F_B$.
We denote them
$\overline{\Cal{L}}_{\tau}$. They are smooth. The respective orbit
which is an open subset of $\overline{\Cal{L}}_{\tau}$ is denoted
$\Cal{L}_{\tau}$.  

\medskip

{\bf 2.6.1. Forgetful morphisms and a family of pointed curves over $\overline{\Cal{L}}_B$.}
Assume that $B\subset B^{\prime}$. Then we have the projection
morphism $\bold{Z}^{B^{\prime}}\to\bold{Z}^B$ which induces the morphism
$f^{B^{\prime},B}:\,N_{B^{\prime}}\to N_B.$ It satisfies the property
stated in the last lines of [Fu], p. 22: for each cone 
$C(\tau^{\prime})\in F_{B^{\prime}}$,
there exists a cone $C(\tau )\in F_{B}$  such that
$f^{B^{\prime},B}(C(\tau^{\prime}))\subset C(\tau ).$ In fact,
$\tau$ is obtained from $\tau^{\prime}$ by deleting
elements of $B^{\prime}\setminus B$ and then deleting the
empty subsets of the resulting partition of $B$.

\smallskip

Therefore, we have a morphism $f^{B^{\prime},B}_*:\,  
\overline{\Cal{L}}_{B^{\prime}}\to\overline{\Cal{L}}_B$
([Fu], p. 23) which we will call {\it forgetful one} (it forgets
elements of $B^{\prime}\setminus B$).

\smallskip

\proclaim{\quad 2.6.2. Proposition} If $B^{\prime}\setminus B$
consists of one element, then the forgetful morphism
$\overline{\Cal{L}}_{B^{\prime}}\to\overline{\Cal{L}}_B$
has a natural structure of a stable $(\{0,\infty\},B)$--pointed curve
of genus zero.
\endproclaim

\smallskip

{\bf Proof.} Let us first study the fibers of the forgetful
morphism. Let $\tau$ be a partition of $B$ of length $l+1$
and $\Cal{L}_{\tau}$ the respective orbit in $\overline{\Cal{L}}_B.$
Its inverse image in $\overline{\Cal{L}}_{B^{\prime}}$ is contained
in the union  $\cup \overline{\Cal{L}}_{\tau^{\prime}}$ where
$\tau^{\prime}$ runs over partitions of $B^{\prime}$ obtained by adding the forgotten point
either to one of the parts $\tau_i$, or inserting it
in between $\tau_i$ and $\tau_{i+1}$, or else
putting it at the very beginning or at the very end as
a separate part.

\smallskip

The inverse image of any point $x\in\Cal{L}_{\tau}$ is acted
upon by the multiplicative group $\bold{G}_m =\roman{Ker}\,(T_{B^{\prime}}\to T_B)$.
This action breaks the fiber into a finite number of orbits
which coincide with the intersections of this fiber
with various $\Cal{L}_{\tau^{\prime}}$ described above.
When $\tau^{\prime}$ is obtained by adding the forgotten
point to one of the parts, this intersection is a torsor
over the kernel, otherwise it is a point. As a result, we get
that the fiber is a chain of $\bold{P}^1$'s, whose components are labeled
by the components of $\tau$ and singular points by the neighboring
pairs of components. 

\smallskip

The forgetful morphism is flat, because locally in toric
coordinates it is described as adjoining a variable and localization.

\smallskip

In order to describe the two white sections of the forgetful
morphism, consider two partitions $(B^{\prime}\setminus B,B)$
and $(B,B^{\prime}\setminus B)$ of $B^{\prime}$ and the respective
closed strata. It is easily seen that the forgetful morphism
restricted to these strata identifies them with $\overline{\Cal{L}}_B$.
We will call them $x_0$ and $x_{\infty}$ respectively.

\smallskip

Finally, to define the $j$--th black section, $j\in B$,
consider the morphism of lattices $s_j:\,N_B\to N_{B^{\prime}}$
which extends a function $\chi$ on $B$ to the function
$s_j(\chi )$ on $B^{\prime}$ taking the value $\chi (j)$ at the forgotten
point. This morphism satisfies the condition of [Fu], p. 22:
each cone $C(\tau )$ from $F_B$ lands in an appropriate cone 
$C(\tau^{\prime})$ from $F_{B^{\prime}}$.
This must be quite clear from the description at the end of 2.5.1: 
$\tau^{\prime}$ is obtained from $\tau$ by adding the forgotten
point to the same part to which $j$ belongs.
Hence we have the induced morphisms $s_{j*}:\,\overline{\Cal{L}}_B\to
\overline{\Cal{L}}_{B^{\prime}}$ which obviously are sections.
Moreover, they do not intersect $x_0$ and $x_{\infty}$, and they
are distributed among the components of the reducible fibers
exactly as expected.

\medskip

\proclaim{ \quad 2.6.3. Theorem} The morphism $\overline{\Cal{L}}_B\to\overline{L}_B$ inducing the family described
in the Proposition 2.6.2 is an isomorphism. 
\endproclaim

\smallskip

This can be proved by induction on $\roman{card}\,B$ with the
help of the more detailed analysis of the forgetful morphism, as above.
We omit the details because they are not instructive.

\smallskip

An important corollary of this Theorem is the existence
of a surjective birational morphism 
$\overline{M}_{0,n+2}\to\overline{L}_n$ corresponding to any choice
of two different labels $i,j$ in $(1,\dots ,n+2)$. In terms of the
of the respective functors, this morphism blows down all the components
of a stable $(n+2)$--labeled curve except for those that belong to
the single path from the component containing the $i$--th point to the one
containing the $j$--th point.

\smallskip

In fact, M.~Kapranov has shown the existence of such a morphism
for $\overline{\Cal{L}}_n$ in place of $\overline{L}_n$ (see [Ka2], p. 102).
He used a different description of $\overline{\Cal{L}}_n$
in terms of the defining polyhedron, which he identified
with the so called permutohedron, the convex hull
of the $\bold{S}_n$--orbit of $(1,2,\dots ,n)$. He has also proved
that $\overline{\Cal{L}}_n$ can be identified with the closure
of the generic orbit of the torus in the space of complete flags
in an $n$--dimensional vector space. 

\medskip

{\bf 2.7. Combinatorial model of $H^*(\overline{\Cal{L}}_B)$.}
We will denote by $[\overline{\Cal{L}}_{\sigma}]_*$
(resp. $[\overline{\Cal{L}}_{\sigma}]^*$) the homology
(resp. the dual cohomology) class of $\overline{\Cal{L}}_{\sigma}$.

\smallskip

The remaining parts of this section (and the Appendix)
are dedicated to the study of linear and non--linear
relations between these classes, in the spirit of
[KM] and [KMK], but with the help of the standard toric techniques. 

\smallskip

Consider a family of pairwise commuting independent variables
$l_{\sigma}$ numbered by 2--partitions of $B$ and introduce
the ring
$$
H^*_B:=\Cal{R}_B/I_B
\eqno(2.13)
$$
where $\Cal{R}_B$ is freely generated by $l_{\sigma}$
(over an arbitrary coefficient ring $k$), and the ideal $I_B$
is generated by the following elements indexed by pairs $i,j\in B$:
$$
r^{(1)}_{ij}:=
 \sum_{\sigma :\,i\sigma j} l_{\sigma}-
\sum_{\tau :\,j\tau i} l_{\tau} ,
\eqno(2.14)
$$
$$
r^{(2)}(\sigma ,\tau ):= l_{\sigma} l_{\tau}\qquad
\roman{if}\ i\sigma j\ \roman{and}\ j\tau i\ \roman{for\ some}\ i,j. 
\eqno(2.15)
$$
\medskip

\proclaim{\quad 2.7.1. Theorem} a) There is a well defined
ring isomorphism $\Cal{R}_B/I_B\to A^*(\overline{\Cal{L}}_B,k)$
such that $l_{\sigma}\,\roman{mod}\,I_B\,\mapsto [\overline{\Cal{L}}_{\sigma}]^*$.
The Chow ring
$A^*(\overline{\Cal{L}}_B,k)$
and the cohomology ring  $H^*(\overline{\Cal{L}}_B,k)$ are canonically isomorphic.

\smallskip

b) The boundary divisors (strata corresponding to 2--partitions) intersect transversally.
\endproclaim

\smallskip

{\bf Proof.} We must check that the ideal of relations between
$2^n-2$ dual classes of the boundary divisors
$[\overline{\Cal{L}}_{\sigma}]^*$ contains and is generated by the following
relations:
$$
R^{(1)}_{ij}:\qquad \sum_{\sigma :\,i\sigma j} [\overline{\Cal{L}}_{\sigma}]^*-
\sum_{\tau :\,j\tau i} [\overline{\Cal{L}}_{\tau}]^* =0.
\eqno(2.16)
$$
If $i\sigma j$ and $j\tau i$, then
$$
R^{(2)}(\sigma ,\tau ):\qquad [\overline{\Cal{L}}_{\sigma}]^*\cdot [\overline{\Cal{L}}_{\tau}]^* =0.
\eqno(2.17)
$$ 
We refer to the Proposition on p. 106 of [Fu] which
gives a system of generators for this ideal for any smooth
proper toric variety (Fulton additionally assumes projectivity
which we did not check, but see [Da], Theorem 10.8 for the general
proper case).

\smallskip

In our notation, these generators look as follows.

\smallskip

To get the complete system of linear relations, we
must choose some elements $m$ in the dual lattice
of $N_B$ spanning this lattice and form the sums
$\sum_{\sigma}m(\chi_{\sigma_1})[\overline{\Cal{L}}_{\sigma}]^*$,
where $\sigma$ runs over all 2--partitions.
In our case, the dual lattice is spanned by the linear
functionals $m_{ij}:\,\chi \mapsto \chi\,(i) -\chi\,(j)$ 
for all pairs $i,j\in B.$
Writing the respective relation,
we get (2.16).

\smallskip

The complete system of nonlinear relations is given by the
monomials $l_{\sigma^{(1)}}\dots l_{\sigma^{(k)}}$ such that
$(C(\sigma^{(1)}),\dots ,C(\sigma^{(k)}))$ do not span a cone
in $F_B$. This means that some pair $(C(\sigma^{(a)}),C(\sigma^{(b)}))$
already does not span a cone, because otherwise the respective
2--partitions would form a good family (cf. 2.4.1). And in view
of Lemma 2.4.2 (iii), we can find $i,j\in B$ such that
$i\sigma^{(a)}j$ and $j\sigma^{(b)}i$. Hence (2.16) and (2.17)
together constitute a generating system of relations.

\smallskip

The remaining statements are true for all smooth complete
toric varieties defined by simplicial fans.

\medskip

{\bf 2.8. Combinatorial structure of the cohomology ring.} In the
remaining part of this section
we fix a finite set $B$ and study $H^*_B$ as an abstract ring.

\smallskip

For an $(N+1)$--partition $\tau$ define the respective {\it good monomial} $m(\tau )$
by the formula
$$
m(\tau )=l_{\sigma^{(1)}}\dots l_{\sigma^{(N)}}\in \Cal{R}_B . 
$$
If $\tau$ is the trivial 1--partition, we put $m(\tau ):=1.$
In view of the Theorem 2.7.1, $m(\tau )$ represents the
cohomology class of $\overline{\Cal{L}}_{\tau}$.

\smallskip

Notice that if we have two good families
of 2--partitions whose union is also good,
then the product of the respective good monomials 
is a good monomial. This defines a partial operation
$*$ on pairs of partitions
$$
m(\tau^{(1)})\,m(\tau^{(2)})=m(\tau^{(1)}*\tau^{(2)}).
$$
\smallskip

\proclaim{\quad 2.8.1. Proposition} Good monomials and $I_B$
span $\Cal{R}_B.$ Therefore, images of good
monomials span $H^*_B.$
\endproclaim

\smallskip

{\bf Proof}. We make induction on the degree. In degrees zero and one
the statement is clear because $l_{\sigma}$ are good.
If it is proved in degree $N$, it suffices to check
that for any 2--partition $\sigma$ and any
nontrivial partition $\tau$, $l_{\sigma}m(\tau )$ is a linear
combination of good monomials modulo $I_B.$ We will consider
the three cases of Lemma 2.4.2 in turn.

\smallskip

(i) {\it $\sigma$ breaks $\tau$ between $\tau_{a}$ and $\tau_{a+1}$}.

\smallskip

This means that $l_{\sigma}$ divides $m(\tau )$.

\smallskip

Choose $i\in \tau_a, j\in \tau_{a+1}.$ In view
of (2.14), we have
$$
\left(\sum_{\rho :\,i\rho j} l_{\rho}-
\sum_{\rho :\,j\rho i} l_{\rho}\right) m(\tau )\in I_B.
\eqno(2.18)
$$
But if $j\rho i$, then $l_{\rho}m(\tau )\in I_B$ because of (2.15).
 Among the terms with
$i\rho j$ there is one $\l_{\sigma}.$ For all other
$\rho$'s, $l_{\rho}$ cannot divide $m(\tau )$ since other
divisors put $i$ and $j$ in the same part of the respective
partition. Therefore, $l_{\rho}m(\tau )$ either belongs
to $I_B$, or is good. So finally (2.18) allows us to
express $l_{\sigma}m(\tau )$ as a sum of good monomials
and an element of $I_B:$
$$
l_{\sigma}m(\tau )= -\sum_{\rho\ne\sigma , \,i\rho j} m(\rho *\tau )
\ \roman{mod}\, I_B
$$
where the terms for which $\rho *\tau$ is not defined must be interpreted as zero. More precisely, there are two types of non--vanishing terms.
One corresponds to all 2--partitions $\alpha$ of $\tau_a$
such that $i\in\tau_{a1}$ which we will write as $i\alpha$.
Another corresponds to 2--partitions $\beta$ of $\tau_{a+1}$
with $j$ belonging to the second part, $\beta j$:
$$
l_{\sigma}m(\tau )= -\sum_{\alpha :\,i\alpha}m(\tau (\alpha))
-\sum_{\beta :\,\beta j}m(\tau (\beta)) \ \roman{mod}\, I_B .
\eqno(2.19)
$$
Notice that there are several ways to write the
right hand side, depending on the choice of $i,\,j.$
Hence good monomials are not linearly independent modulo
$I_B.$

\medskip

(ii) {\it $\sigma$ breaks $\tau$ at $\tau_{a}$}.

\smallskip

According
to the analysis above, this means that 
$$
l_{\sigma}m(\tau )=m(\sigma *\tau )=m(\tau (\alpha))
\eqno(2.20)
$$
for an appropriate partition $\alpha$ of $\tau_a$.

\medskip

(iii) {\it $\sigma$ does not break $\tau$}.

\smallskip

In this case, let $(\tau_b,\tau_{b+1})$ be a bad pair
for $\sigma$. Then from (2.12) it follows that 
there exist $i,j\in B$ such that
$i\sigma j$ and $j \sigma^{(a)}i$. Hence $l_{\sigma}m(\tau )$
is divisible by $r^{(2)}(\sigma ,\sigma^{(a)})$ and 
$$
l_{\sigma}m(\tau ) = 0\ \roman{mod}\,I_B.
$$
\medskip

{\bf 2.8.2. Linear combinations of good monomials
belonging to $I_B$.} Let $\tau =(\tau_1,\dots ,\tau_{N})$ be a partition of $B$. Choose $a\le N$ such that $|\tau_a|\ge 2$,
and two elements $i,j\in \tau_a,\ i\ne j.$ For any ordered
2--partition $\alpha = (\tau_{a1},\tau_{a2})$ of $\tau_a$, denote by $\tau (\alpha )$
the induced $N+1$--partition of $B$ as above:
$$
(\tau_1,\dots ,\tau_{a-1},\tau_{a1},\tau_{a2},\tau_{a+1},\dots ,\tau_N) .
$$
Finally, put
$$
r^{(1)}_{ij}(\tau , a):=\sum_{\alpha :\,i\alpha j} m(\tau (\alpha )) -
\sum_{\alpha :\,j\alpha i} m(\tau (\alpha )).
\eqno(2.21)
$$
Choosing for $\tau$ the trivial 1--partition, we get (2.14)
so that these elements span the intersection
of $I_B$ with the space of good monomials of degree one.

\smallskip

\smallskip

Generally, all $r^{(1)}_{ij}(\tau , a)$ belong to $I_B.$
In fact, keeping the notations above, consider
$$
 r_{ij}^{(1)}m(\tau )= \left( \sum_{\rho :\,i\rho j} l_{\rho}-
\sum_{\rho :\,j\rho i} l_{\rho}\right) m(\tau ) \in I_B .
\eqno(2.22)
$$
Arguing as above, we see that the summand corresponding
to $\rho$ in (2.18) either belongs to $I_B$, or is a good
monomial, and the latter happens exactly for those
partitions $\rho$ that are of the type $\tau (\alpha )$
with either $i\alpha j$, or $j\alpha i.$
Hence (2.21) lies in $I_B.$ This proves our claim.

\medskip

\proclaim{\quad 2.9. Theorem} Elements (2.21) span
the intersection of $I_B$ with the space generated by good monomials.
\endproclaim

\smallskip
 
{\bf Proof.} Define the linear space $H_{*B}$ generated
by the symbols $\mu (\tau )$ for all partitions of $B$ as above
which satisfy analogs of the linear relations (2.21): for all
$(\tau,\tau_a,i,j)$ as above we have
$$
\sum_{\alpha :\,i\alpha j} \mu (\tau (\alpha )) -
\sum_{\alpha :\,j\alpha i} \mu (\tau (\alpha )) =0.
\eqno(2.23)
$$
\proclaim{\quad 2.9.1. Technical Lemma} There exists an (obviously unique) structure of the
$H^*_B$--module on $H_{*B}$ with the following multiplication table.

\smallskip

(i) If $\sigma$ breaks $\tau$ between $\tau_a$ and $\tau_{a+1}$,
then  for any choice of $i\in \tau_a,\,j\in\tau_{a+1}$ 
$$
l_{\sigma}\mu (\tau )= -\sum_{\alpha :\,i\alpha}\mu (\tau (\alpha))
-\sum_{\beta :\,\beta j}\mu (\tau(\beta )).
\eqno(2.24)
$$
(cf. (2.19)).
\smallskip

(ii) If $\sigma$ breaks $\tau$ at $\tau_a$, then
$$
l_{\sigma}\mu (\tau )=\mu (\sigma *\tau ).
\eqno(2.25)
$$ 
(cf. (2.20)).
\smallskip

(iii) If $\sigma$ does not break $\tau$, then
$$
l_{\sigma}\mu (\tau ) =0.
\eqno(2.26)
$$
\endproclaim

\medskip

Our proof of the Technical Lemma consists in the direct verification
that the prescriptions (2.24)--(2.23) are compatible with
all relations that we have postulated. 
Unfortunately, such strategy requires the painstaking case--by--case treatment
of a long list of combinatorially distinct situations,
and we relegate it to the Appendix.

\smallskip

{\bf 2.9.2. Deduction of Theorem 2.9 from the Technical Lemma.} 
Since elements (2.21) belong to $I_B,$ there exists a surjective linear map
$s:\,H_{*B}\to H^*_B$, $\mu (\tau )\mapsto m(\tau ).$ Now denote
by $\bold{1}$ the element $\mu (\varepsilon )$ where $\varepsilon$ is
the 1--partition. Then $t:\,m(\sigma )\mapsto m(\sigma )\bold{1}$
is a linear map $H^*_B\to H_{*B}$. From (2.25) one easily
deduces that $m(\tau )\bold{1} =\mu (\tau )$ so that
$s$ and $t$ are mutually inverse. Therefore, (2.22) span
the linear relations between the images of good monomials in
$H^*_B.$

\medskip

According to the Theorem 2.4.1, $H_{*B}$, together with its structure
of $H^*_B$--module, is a combinatorial model of the
homology module $H_*(\overline{\Cal{L}}_B,k)$. The generators
$\mu (\tau )$ correspond to $[\overline{\Cal{L}}_{\tau}]_*$.

\newpage

\centerline{\bf \S 3. Pencils of flat connections}

\smallskip

\centerline{\bf and the Commutativity Equations}

\medskip

{\bf 3.1. Notation.} Let $M$ be a (super)manifold over a field
$k$ of characteristic zero in one of the
standard categories (smooth, complex analytic, schemes, formal $\dots$).
We use the conventions spelled out in [Ma], I.1. In particular,
differentials in the de Rham complex $(\Omega^*_M,d)$ and connections are odd. 
This determines our sign rules; parity of an object $x$ is denoted
$\widetilde{x}.$

\smallskip

Let $\Cal{F}$ be a locally free sheaf (of sections of a vector bundle)
on $M,$ $\nabla_0$ a connection on $\Cal{F}$, that is an odd
$k$--linear operator $\Cal{F}\to\Omega^1_M\otimes\Cal{F}$ satisfying the Leibniz identity
$$
\nabla_0(\varphi f)=d\varphi\otimes f+(-1)^{\widetilde{\varphi}}
\varphi\,\nabla_0f, \ \varphi\in \Cal{O}_M,\,f\in\Cal{F}.
\eqno(3.1) 
$$
This operator extends to a unique operator on the $\Omega_M^*$--module
$\Omega_M^*\otimes \Cal{F}$ denoted again $\nabla_0$ and satisfying
the same identity (3.1) for any $\varphi\in\Omega_M$. Any other connection
differential $\nabla$ restricted to $\Cal{F}$ has the form
$\nabla_0+\Cal{A}$ where $\Cal{A}:\,\Cal{F}\to\Omega^1_M\otimes \Cal{F}$
is an odd $\Cal{O}_M$--linear operator: $\Cal{A} (\varphi f)=
(-1)^{\widetilde{\varphi}}\varphi\Cal{A}(f).$ Any connection
naturally extends to the whole tensor algebra generated by $\Cal{F},$
in particular, to $\Cal{E}nd\,\Cal{F}.$

\smallskip

The connection
$\nabla_0$ is called flat, iff $\nabla_0^2=0.$ {\it A pencil
of flat connections} is a line in the space of connections
$\nabla_{\lambda}:=\nabla_0 +\lambda\Cal{A}$ such that
$\nabla_{\lambda}^2=0$ ($\lambda$ is an even parameter).
In the smooth, analytic or formal category, $\nabla_0$
is flat iff $\Cal{F}$ locally admits a basis of flat
sections $f, \nabla_0f=0$.

\medskip

\proclaim{\quad 3.2. Proposition} $\nabla_0 +\lambda\Cal{A}$
is a pencil of flat connections iff the following two conditions
are satisfied:

\smallskip

(i) Everywhere locally on $M$, we have 
$$
\Cal{A}=\nabla_0\Cal{B}
\eqno(3.2)
$$
for some $\Cal{B}\in\Cal{E}nd\,\Cal{F}$.

\smallskip

(ii) Such an operator $\Cal{B}$ satisfies the quadratic differential
equation
$$
\nabla_0\Cal{B}\wedge\nabla_0\Cal{B} =0.
\eqno(3.3)
$$
\endproclaim

{\bf Proof.} Calculating the coefficient of $\lambda$
in $\nabla_{\lambda}^2=0$ we get $\nabla_0\Cal{A}=0$. But the complex
$\Omega_M^*\otimes\Cal{F}$ is the resolution of the
sheaf of flat sections $\roman{Ker}\,\nabla_0\subset\Cal{F}.$
This furnishes (i); (ii) means the vanishing of the
coefficient at $\lambda^2.$

\medskip

{\bf 3.2.1. Remarks.} a) Write $\Cal{B}$ as a matrix in a 
basis of $\nabla_0$--flat sections of $\Cal{F}$, whose entries
are local functions on $M.$ Then (3.3) becomes
$$
d\Cal{B}\wedge d\Cal{B}=0.
\eqno(3.4)
$$
These equations written in local coordinates $(t^i)$ on $M$ were called
``$t$--part of the $t$--$t^*$ equations'' by S.~Cecotti
and C.~Vafa. A.~Losev in [Lo1] suggested to call them
``the Commutativity Equations''.

\smallskip

b) If $\nabla_0\varphi_0=0,$ then
$$
(\nabla_0+\lambda\nabla_0\Cal{B})(e^{-\lambda\Cal{B}}\varphi_0)=0.
$$

\medskip

{\bf 3.2.2. Pencils of flat connections related to
Frobenius manifolds.} Any solution to the Associativity Equations produces
a pencil of flat connections. 

\smallskip

To explain this we will use the
geometric language due to B.~Dubro\-vin (and the notation of [Ma], I.1.5).
Consider a Frobenius manifold $(M,g, \circ )$ where  
$$
\circ :\,\Cal{T}_M\otimes_{\Cal{O}_M}\Cal{T}_M\to \Cal{T}_M
$$
is a (super)commutative associative multiplication on the tangent
sheaf satisfying the potentiality condition, and 
$g$ is an invariant flat metric (no positivity condition is assumed,
only symmetry and non--degeneracy). 
Denote by $\nabla_0$ the Levi--Civita connection of $g$.
Finally, denote by $\Cal{A}$ the operator obtained from
the Frobenius multiplication in $\Cal{T}_M$ ([Ma], I.1.4). In other words,
consider the pencil of connections on $\Cal{F}=\Cal{T}_M$ whose covariant derivatives are
$$
(\nabla_{0}+\lambda\Cal{A})_X(Y):=\nabla_{0,X}(Y)+\lambda\,X\circ Y\,.
$$
This pencil is flat (see [Ma], Theorem I.1.5, p. 20). In fact, $\Cal{B}$ written in a basis
of $\nabla_0$--flat coordinates and the respective flat vector fields
is simply the matrix of the second derivatives of a local potential $\Phi$
(with one subscript raised). This is the first structure connection
of $M$.

\smallskip

This pencil admits an infinite dimensional deformation:
one should take the canonical extension of the potential
to the large phase space and consider the coordinates
with gravitational descendants as parameters of the
deformation.

\smallskip

Another family of flat connections, this time on the {\it cotangent}
sheaf of a Frobenius manifold $M$ admitting an Euler vector field $E$
(see [Ma], pp. 23--24), is defined as follows. Denote
the scalar product on vector fields $\check{g}_{\lambda}(X,Y)
:=g((E-\lambda )^{-1}\circ X,Y).$ The inverse form induces a pencil
of flat metrics on the cotangent sheaf, whose Levi--Civita
connections however do not form a pencil of flat connections 
in our sense (see [Du1], Appendix D,
and [Du3] for a general discussion of such setup). This is the second structure
connection of $M$.

\medskip

{\bf 3.2.3. Flat coordinates and gravitational descendants.} 
One can show that 1--forms on $M$ flat with respect to the dualized first
structure connection are closed and therefore locally exact.
Their integrals are called deformed flat coordinates.
In [Du2], Example 2.3 and Theorem 2.2, B.~Dubrovin gives explicit
formal series in $\lambda$ ($z$ in his notation) for suitably normalized  
deformed flat coordinates. Coefficients of these series involve
some correlators with gravitational descendants, namely
those for which the non--trivial operators $\tau_p$
are applied only at one point. In [KM2] and [Ma], VI.7.2, p.278, 
it was shown
that two--point correlators of this kind determine a linear
operator in the large phase space which transforms
the modified correlators with descendants into non--modified ones (in any genus).
This is important because apriori only modified correlators are defined
for an arbitrary Cohomological Field Theory in the sense of [KM1],
which is not necessarily quantum cohomology of a manifold.

\medskip

{\bf 3.2.4. Pencils of flat connections in a global setting.}
Pencils of flat connections appear also in the context
of Simpson's non--abelian Hodge theory. Briefly, consider
a smooth projective manifold $M$ over $\bold{C}$. One can define
two moduli spaces, $Mod_1$ and $Mod_2$. The first one classifies
flat connections (on variable vector bundles $\Cal{F}$
with vanishing rational Chern classes) with semisimple
Zariski closure of the monodromy group. The second one 
classifies semistable Higgs pairs $(\Cal{F}, \Cal{A})$ where $\Cal{A}$
is an operator as in 3.1, satisfying only the condition
$\Cal{A}\wedge\Cal{A}=0.$ (In fact, one should only consider
smooth points of the respective moduli spaces). N.~Hitchin,
C.~Simpson, Fujiki et al. established that $Mod_1$ and $Mod_2$
are canonically isomorphic as $C^{\infty}$--manifolds,
but their complex structures $I$, $J$ are different,
and together with $K=IJ$ produce a hypercomplex manifold.

\smallskip

P. Deligne has shown that the respective twistor space is
precisely the moduli space of the pencils of flat
connections on $M$ (where the Higgs complex structure
corresponds to the point $\lambda =\infty$ in our notation).

\smallskip

For details, see [Si].

\medskip

{\bf 3.3. Formal solutions to the Commutativity Equations
and the homology of $\overline{L}_n.$} In [KM1] and [KMK]
it was shown that formal solutions to the Associativity Equations
are cyclic algebras over the cyclic genus zero homology modular
operad $(H_*(\overline{M}_{0,n+1}))$ (see also [Ma], III.4).
The main goal of this section is to show the similar
role of the homology of the spaces $\overline{L}_n$
in the theory of the Commutativity Equations.
This was discovered and discussed on a physical level in [Lo1], [Lo2].
Here we supply precise mathematical statements with proofs.

\smallskip

Unlike the case
of the Associativity Equations, we will have to deal here
with modules over an algebra (depending explicitly
on the base space) rather than with algebras over an operad.
The main ingredient of the construction is
the direct sum of the homology spaces
of all $\overline{L}_n$ endowed with the multiplication
coming from the boundary morphisms. We work
with the combinatorial models of these spaces
defined in 2.9.1.  

\smallskip

We start with some preparations. Let $V=\oplus_{n=1}^{\infty}V_n$
be a graded associative $k$--algebra (without identity) in
the category of vector $k$--superspaces over a field $k$.
We will call it  {\it an $\bold{S}$--algebra}, if for each $n$,
an action of the symmetric group $\bold{S}_n$ on $V_n$
is given such that the multiplication map $V_m\otimes V_n\to V_{m+n}$
is compatible with the action of $\bold{S}_m\times\bold{S}_n$
embedded in an obvious way into $\bold{S}_{m+n}.$

\smallskip

If $V$ is an $\bold{S}$--algebra, then the sum of subspaces
$J_n$ spanned by $(1-s)v,\,s\in\bold{S}_n,\,v\in V_n,$
is a double--sided ideal in $V.$ Hence the sum of
the coinvariant spaces $V_{\bold{S}_n} =V_n/J_n$ is a graded ring which
we denote $V_{\bold{S}}$.

\smallskip

If $V$, $W$ are two $\bold{S}$--algebras, then the diagonal
part of their tensor product $\oplus_{n=1}^{\infty}V_n\otimes W_n$
is an  $\bold{S}$--algebra as well.

\smallskip

Let $T$ be a vector superspace (below always assumed 
finite--dimensional). Its tensor algebra
(without the rank zero part) is an  
$\bold{S}$--algebra.

\smallskip

As a less trivial example, consider $H_*:=\oplus_{n=1}^{\infty}
H_{*n}$ where we write $H_{*n}$ for $H_{*\{1,\dots ,n\}}$.
The multiplication law is given by what
becomes the boundary morphisms in the geometric setting:
if $\tau^{(1)}$ (resp. $\tau^{(2)}$)
is a partition of $\{1,\dots ,m\}$ (resp. of $\{1,\dots ,n\}$),
then
$$
\mu (\tau^{(1)})\mu (\tau^{(2)}) = \mu (\tau^{(1)}\cup \tau^{(2)})
\eqno(3.5)
$$
where the concatenated partition of 
$\{1,\dots ,m,\, m+1,\dots ,m+n\}$ is defined in an obvious way,
shifting all the components of $\tau^{(2)}$ by $m$.

\smallskip

Our main protagonist is the algebra of coinvariants
of the diagonal tensor product of these examples:
$$
H_{*}T:= \left(\oplus_{n=1}^{\infty} H_{*n}\otimes T^{\otimes n}
\right)_{\bold{S}} .
\eqno(3.6)
$$
We now fix $T$ and another  vector superspace $F$ and assume that the
ground field $k$ has characteristic zero.

\smallskip

\proclaim{\quad 3.3.1. Theorem} There is a natural
bijection between the set of representations of $H_{*}T$
in $F$ and the set of pencils of flat connections
on the trivial bundle with fiber $F$ on the formal
completion of $T$ at the origin.
\endproclaim

\smallskip

This bijection will be precisely defined and discussed
below: see Proposition 3.6.1. Before passing to this definition and the proof
of the Theorem, we will give a down--to--earth coordinate--dependent
description of the representations of $H_*T$.

\medskip

{\bf 3.4. Matrix correlators.} Fix $T$ and choose
its parity homogeneous basis $(\Delta_a\,|\,a\in I)$
where $I$ is a finite set of indices.

\smallskip

For any $n\ge 1$, the space $H_{*n}\otimes T^{\otimes n}$ is spanned
by the elements 
$$
\mu (\tau^{(n)})\otimes \Delta_{a_1}\otimes\dots\otimes\Delta_{a_n}
\eqno(3.7)
$$ 
where $\tau^{(n)}$ runs over all partitions of $\{1,\dots ,n\}$
whereas $(a_1,\dots ,a_n)$ runs over all maps $\{1,\dots ,n\}\mapsto I:\,
i\to a_i.$

\smallskip

In view of the Theorem 2.9, all linear relations between these
elements are spanned by the following ones: choose $(a_1,\dots ,a_n)$
and $(\tau^{(n)}, \tau^{(n)}_r,\,i\ne j\in \tau^{(n)}_r)$, then
$$
\sum_{\alpha :\,i\alpha j} \mu (\tau^{(n)} (\alpha ))
\otimes \Delta_{a_1}\otimes\dots\otimes\Delta_{a_n} -
\sum_{\alpha :\,j\alpha i} \mu (\tau^{(n)} (\alpha ))
\otimes \Delta_{a_1}\otimes\dots\otimes\Delta_{a_n} =0
\eqno(3.8)
$$
where the summation is taken over all 2--partitions $\alpha$
of $\tau^{(n)}_r$ separating $i$ and $j$.

\smallskip

The action of a permutation $i\mapsto s(i)$ on (3.7) is defined by
$$
s\left(\mu (\tau^{(n)})\otimes \Delta_{a_1}\otimes\dots\otimes\Delta_{a_n}
\right)
= 
\varepsilon (s, (a_i))\, \mu (s(\tau^{(n)}))\otimes \Delta_{a_{s(1)}}\otimes\dots\otimes\Delta_{a_{s(n)}} .
\eqno(3.9)
$$
Here $\varepsilon (s, (a_i))=\pm 1$ is the sign of the permutation
induced by $s$ on the subfamily of odd $\Delta_{a_i}$'s,
and $s(\tau^{(n)})$ is defined as follows:
$$
s(i)\in s(\tau^{(n)})_r\quad \roman{iff}\quad i\in \tau^{(n)}_r .
\eqno(3.10)
$$
Finally, the multiplication rule between the generators
in the diagonal tensor product is given by:
$$
\mu (\tau^{(m)})\otimes \Delta_{a_1}\otimes\dots\otimes\Delta_{a_m}\cdot\,
\mu (\tau^{(n)})\otimes \Delta_{b_1}\otimes\dots\otimes\Delta_{b_n}
$$
$$
=
\mu (\tau^{(m)}\cup\tau^{(n)})\otimes \Delta_{a_1}\otimes\dots\otimes\Delta_{a_m}\otimes
\Delta_{b_1}\otimes\dots\otimes\Delta_{b_n} .
\eqno(3.11)
$$
Any linear representation $K:\,H_*T\to \roman{End}\,F$ can
be described as a linear representation of the diagonal tensor
product satisfying additional symmetry restrictions.
To spell it out explicitly, we define {\it the matrix correlators}
of $K$ as the following family of endomorphisms of $F$:
$$
\tau^{(n)}\langle \Delta_{a_1}\dots\Delta_{a_n}\rangle :=
K(\mu (\tau^{(n)})\otimes \Delta_{a_1}\otimes\dots\otimes\Delta_{a_n}) .
\eqno(3.12)
$$
\smallskip

\proclaim{\quad 3.4.1. Claim} Matrix correlators of any representation
satisfy the following relations:

\smallskip

(i)\ $\bold{S}_n$--symmetry:
$$
s^{-1}(\tau^{(n)})\langle \Delta_{a_1}\dots\Delta_{a_n}\rangle =
\varepsilon (s, (a_i))\, \tau^{(n)}\langle \Delta_{a_{s(1)}}\dots\Delta_{a_{s(n)}}\rangle  \,.
\eqno(3.13)
$$

\smallskip

(ii) Factorization:
$$
(\tau^{(m)}\cup \tau^{(n)})\langle \Delta_{a_1}\dots\Delta_{a_m}
\Delta_{b_1}\dots\Delta_{b_n}\rangle
=
\tau^{(m)}\langle \Delta_{a_1}\dots\Delta_{a_m}\rangle \cdot
\tau^{(n)}\langle \Delta_{b_1}\dots\Delta_{b_n}\rangle \,.
\eqno(3.14)
$$
\smallskip

(iii) Linear relations:
$$
\sum_{\alpha :\,i\alpha j} \tau^{(n)} (\alpha )
\langle \Delta_{a_1}\dots\Delta_{a_n}\rangle -
\sum_{\alpha :\,j\alpha i} \tau^{(n)} (\alpha )
\langle \Delta_{a_1}\dots\Delta_{a_n}\rangle =0
\eqno(3.15)
$$

\smallskip

Conversely, any family of elements of $\roman{End}\,F$
defined for all $n,\,(a_1,\dots ,a_n),\,\tau^{(n)}$ and
satisfying (3.13)--(3.15) consists of matrix correlators
of a well defined representation $K:\,H_*T\to\roman{End}\,F$.
\endproclaim

In fact, we obtain (3.13) by applying $K$ to (3.9)
written for $s^{-1}(\tau^{(n)})$ in place of $\tau^{(n)}$,
because $K$, coming from $H_*T$, vanishes on the image of $1-s$.
Moreover, (3.14) means the compatibility with the multiplication of the
generators.
Finally, (3.15) is a necessary and sufficient condition
for the extendability of the system of matrix correlators to
a linear map $K$.

\medskip

Notice that we can replace here $\roman{End}\,F$ by an
arbitrary associative superalgebra over $k$.

\medskip

{\bf 3.5. Top matrix correlators.} Define {\it top matrix
correlators of $K$} as the subfamily of correlators corresponding
to the identical partitions $\varepsilon^{(n)}$ of $\{1,\dots ,n\}$:
$$
\langle \Delta_{a_1}\dots\Delta_{a_n}\rangle :=
\varepsilon^{(n)}\langle \Delta_{a_1}\dots\Delta_{a_n}\rangle\, .
$$
\smallskip

\proclaim{\quad 3.5.1. Proposition} Top matrix correlators satisfy
the following relations:
$$
\langle \Delta_{a_1}\dots\Delta_{a_n}\rangle =
\varepsilon (s, (a_i))\,\langle \Delta_{a_{s(1)}}\dots\Delta_{a_{s(n)}}\rangle  
\eqno(3.16)
$$
and
$$
\sum_{\sigma:\,i\sigma j}\varepsilon (\sigma,(a_k))
\langle\prod_{k\in\sigma_1}\Delta_{a_k}\rangle\cdot
\langle\prod_{k\in\sigma_2}\Delta_{a_k}\rangle - 
\sum_{\sigma:\,j\sigma i}\varepsilon (\sigma,(a_k))
\langle\prod_{k\in\sigma_1}\Delta_{a_k}\rangle\cdot
\langle\prod_{k\in\sigma_2}\Delta_{a_k}\rangle =0\,.
\eqno(3.17)
$$
Here $\sigma$ runs over 2--partitions of $\{1,\dots ,n\}$.
We choose additionally an arbitrary ordering of both parts $\sigma_1,
\sigma_2$ determining the ordering of $\Delta$'s in the
angular brackets, and compensate this choice by the
$\pm 1$--factor $\varepsilon (\sigma,(a_k))$.

\medskip

Conversely, any family of elements 
$\langle\Delta_{a_1}\dots\Delta_{a_n}\rangle\in \roman{End}\,F$ defined for all $n$ and
$(a_1,\dots ,a_n)$ and satisfying (3.16), (3.17)
is the family of top matrix correlators of a well defined
representation  $K:\,H_*T\to \roman{End}\,F$.
\endproclaim 

\smallskip

{\bf Proof.} Clearly, (3.16) is a particular case of (3.13).
To get (3.17), we apply (3.15) to the identical partition
$\tau^{(n)}=\varepsilon^{(n)}$ and then replace each term 
by the double product of top correlators using (3.14).

\smallskip

Conversely, assume that we are given $\langle \Delta_{a_1}\dots\Delta_{a_n}\rangle$ satisfying (3.16) and (3.17).
There is a unique way to extend this system to a family
of elements $\tau^{(n)}
\langle \Delta_{a_1}\dots\Delta_{a_n}\rangle$ defined for
all $N$--partitions $\tau^{(n)}$ and satisfying the
factorization property (3.14) and at least a part of the symmetry
relations (3.13):
$$
\tau^{(n)}
\langle \Delta_{a_1}\dots\Delta_{a_n}\rangle :=
\varepsilon (\tau^{(n)},(a_k))
\prod_{r=1}^{r=N}\langle\prod_{k\in\tau^{(n)}_r}\Delta_{a_k}\rangle \, .
\eqno(3.18)
$$
Here, as in (3.17), we choose  arbitrary orderings
of each $\tau^{(n)}_r$ and compensate this by the appropriate sign
so that the result does not depend on the choices made. 
All the relations (3.13) become automatically satisfied 
with this definition. In fact, the left hand side
of (3.13) puts into $s^{-1}(\tau^{(n)})_r$
those $i$ for which $s(i)\in \tau^{(n)}_r$ (see (3.10))
so that the expression of both sides of (3.13) through the top correlators
consists of the same groups taken in the same order.
The coincidence of the signs is left to the reader.

\smallskip

It remains to check that (3.18) satisfy the linear relations
(3.15). Recall now that to write a concrete relation (3.15) down
we choose $\tau^{(n)},\, r$, $i,j\in\tau_r^{(n)}$ and $(a_1,\dots ,a_n)$ 
and then sum
over 2--partitions $\alpha$ of $\tau_r^{(n)}$. Hence replacing
each term of the left hand side of (3.15) by the prescriptions (3.18) we get
$$
\prod_{p=1}^{r-1}\langle\prod_{k\in\tau^{(n)}_p}\Delta_{a_k}\rangle\cdot
\left( \sum_{\alpha:\,i\alpha j}\pm\langle\prod_{k\in\alpha_1} \Delta_{a_k}\rangle\cdot
\langle\prod_{k\in\alpha_2}\Delta_{a_k}\rangle -
\sum_{\alpha:\,j\alpha i}\pm\langle\prod_{k\in\alpha_1} \Delta_{a_k}\rangle\cdot
\langle\prod_{k\in\alpha_2}\Delta_{a_k}\rangle
   \right)
$$
$$
\cdot
\prod_{q=r+1}^{N}\langle\prod_{k\in\tau^{(n)}_q}\Delta_{a_k}\rangle \,.
$$
This expression vanishes because its middle term is an instance of (3.17).

\medskip

{\bf 3.6. Precise statement and proof of the Theorem 3.3.1.} 
Assume that we are given a representation $K:\,H_*T\to\roman{End}\,F.$
We will produce from it a formal solution of the Commutativity Equations
using only its top correlators. Let $(x^a)$ be the basis of
formal coordinates on $T$ dual to $(\Delta_a)$. Put
$$
\Cal{B} =
\sum_{n=1}^{\infty}\sum_{(a_1,\dots ,a_n)}
\frac{x^{a_n}\dots x^{a_1}}{n!}\,
\langle \Delta_{a_1} \dots \Delta_{a_n}\rangle \in k[[x]]\otimes \roman{End}\,F.
\eqno(3.19)
$$

\proclaim{\quad 3.6.1. Proposition} a) We have
$$
d\Cal{B}\wedge d\Cal{B}=0.
\eqno(3.20)
$$

\smallskip

b) Conversely, let $\Delta (a_1,\dots ,a_n)\in\roman{End}\,F$ be a
family of linear operators defined for all $n\ge 1$ and all
maps $\{1,\dots ,n\}\to I:\,i\mapsto a_i$. Assume that the parity
of $\Delta (a_1,\dots ,a_n)$ coincides with the sum of the
parities of $\Delta_{a_i}$ and that for any $s\in\bold{S}_n$
$$
\Delta (a_{s(1)},\dots ,a_{s(n)})=\varepsilon (s,(a_i))\,\Delta (a_1,\dots ,a_n)\,.
$$
Finally, assume that the formal series 
$$
\Cal{B} =
\sum_{n=1}^{\infty}\sum_{(a_1,\dots ,a_n)}
\frac{x^{a_n}\dots x^{a_1}}{n!}\,
\Delta ({a_1}, \dots ,{a_n}) \in k[[x]]\otimes \roman{End}\,F
\eqno(3.21)
$$
satisfies the equations (3.20). Then there exists a well defined
representation $K:\,H_*T\to \roman{End}\,F$ such that
$\Delta ({a_1}, \dots ,{a_n})$ are the top correlators
$\langle \Delta_{a_1} \dots \Delta_{a_n}\rangle$ of this
representation.
\endproclaim

\smallskip

Notice that any even element of $k[[x]]\otimes\roman{End}\,F$
without constant term can be uniquely written in the
form (3.21).

\smallskip

{\bf Proof.} Clearly, the equations $d\Cal{B}\wedge d\Cal{B}=0$
written for the series (3.21) are equivalent to a family of
bilinear relations between the symmetric matrix--valued tensors
$\Delta ({a_1} \dots {a_n})$. In view of the Proposition 3.5.1,
it remains to check only that this family of relations
is equivalent to the family (3.17). This is a straightforward
exercise.

\newpage

\centerline{\bf \S 4. Stacks $\overline{L}_{g;A,B}$
and the extended modular operad}

\medskip

{\bf 4.1. Introduction.} The basic topological operad $(\overline{M}_{0;n+1},
n\ge 2)$ of Quantum Cohomology lacks the $n=1$ term which
is usually formally defined as a point. We argued elsewhere 
(cf. [MaZ], sec. 7 and [Ma], VI.7.6) that
it would be very desirable to find a non--trivial DM--stack
which could play the role of $\overline{M}_{0;2}$.
There are several tests that such an object should pass:

\smallskip

a) It must be a semigroup (because for any operad $\Cal{P}$,
the operadic multiplication makes a semigroup of $\Cal{P}(1)$).

\smallskip

b) It must be a part of an extended genus zero operad,
say, $(\widetilde{L}_{0;n+1},
n\ge 1)$ geometrically related to $(\overline{M}_{0;n+1},
n\ge 2)$ in such a way that the theory of Gromov--Witten
invariants with gravitational descendants could be
formulated in this new context. In particular, it must
geometrically explain two--point correlators
with gravitational descendants.

\smallskip

c) In turn, the extended genus zero operad must be a part
of an extended modular operad containing moduli spaces of
arbitrary genus, in such a way that algebras over
classical modular operads produce extended algebras.

\medskip

In this section we will try to show that the space
$$
\widetilde{L}_{0;2}:=\coprod_{n\ge 1} \overline{L}_n 
\eqno(4.1)
$$
passes at least a part of these tests. (Another candidate
which might be interesting is $\roman{lim\,proj}\,\overline{L}_n$
with respect to the forgetful morphisms). 

\medskip

{\bf 4.2. Semigroup structure.} It is defined
as the union of boundary (clutching) morphisms
$$
b:= (b_{n_1,n_2}):\ \widetilde{L}_{0;2}\times \widetilde{L}_{0;2}\to \widetilde{L}_{0;2}
\eqno(4.2)
$$
where 
$$
b_{n_1,n_2}:\ \overline{L}_{n_1}\times \overline{L}_{n_2}\to
\overline{L}_{n_1+n_2}
$$
glues $x_{\infty}$ of the first curve to $x_0$ of the second curve
and renumbers the black points of the second curve keeping
their order (cf. [MaZ], section 7). This is the structure that
induced our multiplication on $H_*$ in 3.3 above.

\medskip

{\bf 4.3. Extended operads.} In (4.2), only
white points $\{x_0,x_{\infty}\}$ are used to define
the operadic composition whereas the black ones
serve only to stabilize the strings of $\bold{P}^1$'s
which otherwise would be unstable. This is a key observation
for our attempt to define an extended operad.

\smallskip

A natural idea would be to proceed as follows. Denote
by $\overline{M}_{g;A,B}$ the stack of stable $(A,B)$--pointed
curves of genus $g$ (see Definition 1.1). Check that it is
a DM--stack. Put $\widetilde{M}_{g;m+1}:=
\coprod_{n\ge 0} \overline{M}_{g;m+1,n}$
 and define the operadic
compositions via boundary maps, using only white points as above.
(We sometimes write here and below $n$ instead of 
$\{1,\dots ,n\}$).

\smallskip

It seems however that this object is too big for our purposes
and that it must be replaced by a smaller stack which we will
define inductively, by using the Construction 1.3
which we will call here simply the adjoining of a generic black point. The components
of this stack will be defined inductively. 

\smallskip

If $g\ge 2$, $m\ge 0$, we start
with $\overline{M}_{g;m}=\overline{M}_{g;m,\emptyset}$ and
add $n$ generic black points, one in turn. Denote the resulting
stack by $\overline{L}_{g;m,n}.$ 

\smallskip

For $g=1$, one should 
add one more sequence of stacks, corresponding to $m=0$.
Since we want to restrict ourselves to Deligne--Mumford stacks, we start at $\overline{M}_{1;0,1}$ identified with $\overline{M}_{1;1}$
(see 1.2 a)), and add black points to get the sequence
$\overline{L}_{1;0,n}$, $n\ge 1.$ These spaces are needed
to serve as targets for the clutching morphisms gluing
$x_0$ to $x_{\infty}$ on the same curve of genus zero: cf. below.

\smallskip

Finally, for $g=0$ we obtain our series of spaces $\overline{L}_n=\overline{L}_{0;2,n}$,
$n\ge 1$ and moreover $\overline{L}_{0;m,n}$, for all $m\ge 3, n\ge 0$.

\medskip

{\bf 4.3.1. Combinatorial types of fibers.} Let us remind
that combinatorial types of classical (semi)stable curves 
with (only white) points labeled by a finite set $A$ 
are isomorphism classes of graphs, whose vertices are
labeled by ``genera'' $g\ge 0$ and tails are bijectively labeled by elements
of $A$. Stability means that vertices of genus 0 bound $\ge 3$
flags, and vertices of genus 1 bound $\ge 1$ flags.
Graphs can have edges with only one vertex, that is, simple loops.
See [Ma], III.2 for more details.

\smallskip

Starting with such a graph $\Gamma$, or rather with its geometric
realization, we can obtain an infinite series
of graphs, which will turn out to be exactly combinatorial
types of (semi)stable $(A,B)$--pointed curves
that are fibers of the families described above. Namely,
subdivide edges and tails of $\Gamma$ by a finite set of new vertices
of genus zero (on each edge or tail, this set may be empty). If
a tail was subdivided, move the respective label (from $A$)
to the newly emerged tail. Distribute the black
tails labeled by elements of $B$ arbitrarily among the old and the new
vertices. Call the resulting graph {\it a stringy stable
combinatorial type} if it becomes stable after repainting
black tails into white ones. Clearly, new vertices depict strings of
$\bold{P}^1$'s stabilized by black points and eventually
two special points on the end components.

\medskip

\proclaim{\quad 4.3.2. Theorem} a) $\overline{L}_{g;m,n}$ is the
Deligne--Mumford stack classifying 
$(m,n)$--pointed curves of genus $g$ of stringy stable combinatorial types.
It is proper and smooth. 
\smallskip

b) Therefore, one can define boundary morphisms gluing two white
points of two different curves
$$
\overline{L}_{g_1;m_1+1,n_1}\times \overline{L}_{g_2;m_2+1,n_2}
\to \overline{L}_{g_1+g_2;m_1+m_2+1,n_1+n_2} 
$$
and gluing two white
points of the same curve:
$$
\overline{L}_{g;m+1,n}\to \overline{L}_{g+1,m-1;n}
$$
such that the locally finite DM--stacks 
$$
\widetilde{L}_{g,m+1}:=
\coprod_{n\ge 0} \overline{L}_{g;m+1,n}
$$
will form components of a modular operad.
\endproclaim

\medskip

The statement a) can be proved in the same way as the
respective statement 2.2 a).

\medskip

It remains to see whether one can develop an extension
of the Gromov--Witten invariants, preferably with descendants, to this
context. The Remark 3.2.3 seems promising in this respect.

\newpage  
   
\centerline{\bf Appendix. Proof of the Technical Lemma}

\bigskip

We break the proof into several steps whose content is indicated
in the title of the corresponding subsection. An advice for the reader
who might care to check the details: the most daunting task
is to convince oneself that none of the alternatives has been inadvertently omitted.

\medskip

{\bf A.1. The right hand side of (2.24) does not depend on
the choice of $i,\,j$.} 

\smallskip

We must check that
a different choice  leads to the same answer modulo
relations (2.23). We can pass from one choice to another by
consecutively replacing only one element of the pair.
Consider, say, the passage from $(i,j)$ to $(i^{\prime},j)$.
Form the difference of the right hand sides of (2.24)
written for $(i^{\prime},j)$ and for $(i,j)$.

\smallskip
In this difference, the terms corresponding
to the partitions $\beta$ will cancel.
The remaining terms will
correspond to the partitions $\alpha$ of $\tau_a$ 
which separate $i$ and $i^{\prime}$. Their difference will vanish in $H_{*B}$
because of (2.23).

\medskip

{\bf A.2. Multiplications by $l_{\sigma}$ are compatible
with linear relations (2.23) between $\mu (\tau )$.}

\smallskip

Choose and fix one linear relation (2.23), that is, a quadruple
$(\tau ,\tau_a, i,j\in \tau_a)$, $i\ne j$. Choose also
a 2--partition $\sigma$. We want to check that after multiplying
the left hand side of (2.23) by $l_{\sigma}$ according
to the prescriptions (2.23)--(2.26) we will get zero modulo
all relations of the type (2.23). There are several basic cases to consider.

\smallskip

{\it (i) $\sigma$ breaks $\tau$ at $\tau_b$, $b\ne a$.} Then put
$\tau^{\prime}=\sigma *\tau$. After multiplication we will
get again (2.23) written for $\tau^{\prime}$ and one of its
components $\tau_a$.

\smallskip

{\it (ii) $\sigma$ breaks $\tau$ at $\tau_a$.} Let $(\tau_{a1},
\tau_{a2})$ be the induced partition; it is now fixed. We must analyze
$l_{\sigma}\mu (\tau (\alpha))$ for variable
2--partitions $\alpha$ of $\tau_a$ with $i\alpha j$ or
$j\alpha i$.

\smallskip

Those $\alpha$ which do not break $(\tau_{a1},\tau_{a2})$
will contribute zero because of (2.26).

\smallskip

Those $\alpha$ which break $(\tau_{a1},\tau_{a2})$
will produce a 3--partition of $\tau_a$, say
$(\tau_{a11},\tau_{a12},\tau_{a2})$ or else
$(\tau_{a1},\tau_{a21},\tau_{a22})$. Finally, there will be
one $\alpha$ which is induced by $\sigma$ that is,
coincides with $(\tau_{a1},\tau_{a2})$.
We must show that the sum total of the respective terms vanishes.
However, the pattern of cancellation will depend
on the positions of $i$ and $j$. In order to present
the argument more concisely, we will first introduce
the numerotation of all possible positions {\it with
respect to a variable $\alpha$} as follows.
Partitions which break $(\tau_{a1},\tau_{a2})$ with $i\alpha j$:
$$
\roman{(I)}:\ i\in\tau_{a11},\,j\in\tau_{a12}\qquad
\roman{(II)}:\ i\in\tau_{a11},\,j\in\tau_{a2}
$$ 
$$
\roman{(III)}:\ i\in\tau_{a1},\,j\in\tau_{a22}\qquad
\roman{(IV)}:\ i\in\tau_{a21},\,j\in\tau_{a22}
$$
Partitions which break $(\tau_{a1},\tau_{a2})$ 
and satisfy $j\alpha i$ will be denoted similarly, but
with prime. Say, (III)${}^{\prime}$ means (III) with
positions of $i$ and $j$ reversed.

\smallskip

Now we will explain the patterns of cancellation depending on the positions
of $i,j$ with respect to $\sigma$. Recall that this latter data
is fixed and determined by the choices we made at the beginning
of this subsection.

\smallskip

If $i,j\in\tau_{a1}$, the only non--vanishing terms
are of the types (I) and (I)${}^{\prime}$. Their sum over
all $\alpha$ will vanish because of (2.23). Similarly,
if $i,j\in\tau_{a2}$, (IV) and (IV)${}^{\prime}$  will cancel,
and everything else will vanish.

\smallskip

Finally, assume that $i\in\tau_{a1},\,j\in\tau_{a2}$. 
that is,  $\sigma$ separates $i,\,j$. Then we may have 
non--vanishing terms of the types (II) and (III)
and in addition the terms coming from (the partition of
$\tau_a$ induced by) $\sigma$ which must be treated
using the formula (2.24), applied however to
$(\tau_1,\dots ,\tau_{a-1},\tau_{a1},\tau_{a2},\tau_{a+1},
\dots )$ in place of $\tau$. Half of these latter terms 
(with $i\in\tau_{a11}$) will
cancel (II), whereas the other half 
(with $j\in\tau_{a22}$) will cancel (III).

\smallskip

The case $j\in\tau_{a1},\,i\in\tau_{a2}$ is treated similarly.

\medskip

{\it (iii) $\sigma$ breaks $\tau$ between $\tau_b$ and $\tau_{b+1}$.}
In this case $\sigma$ breaks any $\tau (\alpha )$ in (2.23)
between two neighbors as well. A contemplation will convince
the reader that only the cases $b=a-1$ and $b=a$ may present
non--obvious cancellations. Let us treat the first one;
the second one is simpler.

\smallskip

For $\alpha =(\tau_{a1},\tau_{a2})$ we will calculate each term 
$l_{\sigma}\mu (\tau (\alpha ))$
using a formula of the type (2.24), first choosing
some $k\in\tau_{a-1},\,l\in\tau_{a1}$ (in place of $i,j$
of (2.24): these letters are already bound). The choice of
$k$ does not matter, but we will choose $l=i$ if $i\alpha j$,
and $l=j$ if $j\alpha i$. We get then for $i\alpha j$:
$$
l_{\sigma}\mu (\tau (\alpha ))=
l_{\sigma}\mu (\dots \tau_{a-1}\tau_{a1}\tau_{a2}\dots )
$$
$$
=-\sum_{\beta :\,k\in\tau_{a-1,1}}
\mu (\dots \tau_{a-1,1}\tau_{a-1,2}\tau_{a1}\tau_{a2}\dots )
-\sum_{\gamma :\,i\in\tau_{a12}}
\mu (\dots \tau_{a-1}\tau_{a11}\tau_{a12}\tau_{a2}\dots )
\eqno(A.1)
$$
where $\beta$ runs over 2--partitions of $\tau_{a-1}$
and $\gamma$ runs over 2--partitions of $\tau_{a1}.$
Write now a similar expression for $j\alpha i$
(with the choice $l=j$). The second sum in this expression
will term--by--term cancel the second sum in (A.1), because
our choices force $i\in\tau_{a12},\,j\in\tau_{a2}$ in both
cases.

\smallskip

If we sum first over $\alpha$, we will see
that the first two sums cancel modulo relations (2.23)
because our choices imply $i\in\tau_{a1},\,j\in\tau_{a2}$ 
in the first sum of (A.1) and the reverse relation
in the first sum written for $j\alpha i.$

\medskip

{\it (iv) $\sigma$ does not break $\tau$.} In this case
we choose a bad pair $(\tau_b,\tau_{b+1})$ for $\sigma$
and $\tau$ (see Lemma 2.4.2(iii)). One easily sees that
if $a\ne b,\,b+1$, then it remains a bad pair for
$\sigma$ and $\tau (\alpha )$ for any $\alpha$ in (2.23).
Therefore, $\l_{\sigma}$ annihilates all terms of
(2.23) in view of (2.26).

\smallskip

We will show that in the exceptional cases we still can find
a bad pair for
$\sigma$ and $\tau (\alpha )$, but it will depend
on $\alpha =(\tau_{a1},\tau_{a2})$, which does not change the remaining
argument.

\smallskip

Assume that $b=a$ that is, $\tau_a\setminus\sigma_1\ne\emptyset ,
\tau_{a+1}\cap\sigma_1\ne\emptyset$ (see (2.12)). Then $(\tau_{a2},\tau_{a+1})$
is a bad pair for $\sigma$ and $\tau (\alpha )$, unless
$\tau_{a2}\subset\sigma_1$, in which case $\sigma_1$
cannot contain $\tau_{a1}$ so that $(\tau_{a1},\tau_{a2})$
form a bad pair.

\smallskip

Similarly, if $b=a-1$, then $\tau_{a-1},\tau_{a1}$ will be
a bad pair unless $\tau_{a1}\cap\sigma_1=\emptyset$,
in which case $(\tau_{a1},\tau_{a2})$ will be a bad pair.

\medskip

By this time we have checked that multiplications by
$l_{\sigma}$ are well defined linear operators on the space
$H_{*B}$. We will now prove that they pairwise commute
and therefore define an action of $\Cal{R}_B$ upon $H_{*B}$.

\medskip

{\bf A.3. Multiplications by $l_{\sigma}$ pairwise commute.}

\smallskip

We start with fixing $\tau$, $\sigma^{(1)}$ and $\sigma^{(2)}$.
We want to check that
$$
l_{\sigma^{(1)}}(l_{\sigma^{(2)}}\mu (\tau )) =
l_{\sigma^{(2)}}(l_{\sigma^{(1)}}\mu (\tau )) .
$$
We may and will assume that $\sigma^{(1)}\ne\sigma^{(2)}$. 
The following alternatives can occur for $\sigma^{(1)}$ and $\sigma^{(2)}$
separately:

\medskip

(i) {\it $\sigma^{(1)}$ breaks $\tau$ at $\tau_a$.}

(ii) {\it $\sigma^{(1)}$ breaks $\tau$ between $\tau_a$ and $\tau_{a+1}$.}

(iii) {\it $\sigma^{(1)}$ does not break $\tau$.}

\smallskip

(i)${}^{\prime}$ {\it $\sigma^{(2)}$ breaks $\tau$ at $\tau_b$.}

(ii)${}^{\prime}$ {\it $\sigma^{(2)}$ breaks $\tau$ between $\tau_b$ and $\tau_{b+1}$.}

(iii)${}^{\prime}$ {\it $\sigma^{(2)}$ does not break $\tau$.}

\medskip

We will have to consider the combined alternatives
(i)(i)${}^{\prime}$,\, (i)(ii)${}^{\prime}$,\, $\dots$ ,
(iii)(iii)${}^{\prime}$ in turn. The symmetry of
$\sigma^{(1)}$ and $\sigma^{(2)}$ allows us to discard
a few of them.

\bigskip

\centerline{{\it Subcase} (i)(i)${}^{\prime}$}

\medskip

We will first assume that $a\ne b$, say $a<b$. Denote
by $\alpha$ (resp. $\beta$) the partition induced
by $\sigma^{(1)}$ (resp. $\sigma^{(2)}$) on $\tau_a$
(resp. $\tau_b$). Then
$$
l_{\sigma^{(1)}}(l_{\sigma^{(2)}}\mu (\tau ))=
l_{\sigma^{(2)}}(l_{\sigma^{(1)}}\mu (\tau ))=
\mu (\tau (\alpha )(\beta )) =\mu (\tau (\beta )(\alpha )) .
$$ 
Now assume that $a=b$. If $\alpha$ breaks $\beta$, we will
again have the desired equality, because $\alpha *\beta =
=\beta *\alpha$. If  $\alpha$ does not break $\beta$,
the both sides will vanish.

\smallskip 

After having treated this subcase, we add one more remark
which will be used below, in the subsection A.5.
Namely, $\alpha$ does not break $\beta$ exactly
when $\sigma^{(1)}$ does not break $\sigma^{(2)}$.
Therefore, if $l_{\sigma^{(1)}}l_{\sigma^{(2)}}$
is one of the quadratic generators of $I_B$,
then consecutive multiplication by the respective
elements annihilates $\mu (\tau )$.

\bigskip

\centerline{{\it Subcase} (i)(ii)${}^{\prime}$}

\medskip

If $a\ne b,\,b+1$, the modifications induced in $\tau$ by
the two multiplications are made in mutually disjoint places
and therefore commute as above. Consider now the case
$a=b$, the case $a=b+1$ being similar. 

\smallskip

Denote by $(\tau_{a1},\tau_{a2})$ the partition induced by $\sigma$
on $\tau_a.$ Then we have
$$
l_{\sigma^{(1)}}\mu (\tau )=
\mu (\dots \tau_{a-1}\tau_{a1}\tau_{a2}\tau_{a+1}\dots )=\mu (\tau^{\prime}).
$$
Clearly, $\sigma^{(2)}$ breaks $\tau^{\prime}$ between $\tau_{a2}$
and $\tau_{a+1}$ so that, after choosing $i\in\tau_{a2},j\in\tau_{a+1}$
we have 
$$ 
l_{\sigma^{(2)}}(l_{\sigma^{(1)}}\mu (\tau ))=-\sum_{\alpha :\,i\alpha}
\mu (\tau^{\prime}(\alpha ))
-\sum_{\beta :\,\beta j}
\mu (\tau^{\prime}(\beta ))
\eqno(A.2)
$$
where $\alpha$ runs over 2--partitions of $\tau_{a2}$
and $\beta$ runs over 2--partitions of $\tau_{a+1}$.

\smallskip

On the other hand, with the same choice of $i,\,j$ we have
$$
l_{\sigma^{(2)}}\mu (\tau )=-\sum_{\gamma :\,i\gamma}
\mu (\tau (\gamma ))
-\sum_{\beta :\,\beta j}
\mu (\tau (\beta ))
\eqno(A.3)
$$
where $\gamma$ runs over 2--partitions of $\tau_{a}$
and $\beta$ runs over 2--partitions of $\tau_{a+1}$.
After multiplication of (A.3) by $l_{\sigma^{(1)}}$,
the second sum in (A.3) will become the second sum of (A.2).
In the first sum, only partitions $\delta$ 
breaking $(\tau_{a1},\tau_{a2})$ will survive,
and they will produce exactly the first sum in (A.2).

\bigskip

\centerline{{\it Subcase} (i)(iii)${}^{\prime}$}

\medskip

Here $\sigma^{(1)}$ breaks $\tau$ at $\tau_a$, and there exists
a bad pair $(\tau_b,\tau_{b+1})$ for $\sigma^{(2)}$
and $\tau$. Since $l_{\sigma^{(2)}}\mu (\tau )=0$,
it remains to check that $l_{\sigma^{(2)}}(l_{\sigma^{(1)}}\mu (\tau ))=0$.
But $l_{\sigma^{(1)}}\mu (\tau ) =\mu (\tau^{\prime})$ as
in the previous subcase. So it remains to find
a bad pair for $\sigma^{(2)}$ and $\tau^{\prime}$.   

\smallskip

If $a\ne b, b+1$, then $(\tau_b,\tau_{b+1})$ will will be such a bad pair

\smallskip

If $a=b$, denote by $(\tau_{a1},\tau_{a2})$ the partition of $\tau_a$
induced by $\sigma^{(1)}$. If $\tau_{a2}$ is not contained
in $\sigma^{(2)}$, $(\tau_{a2},\tau_{a+1})$ will form a bad pair.
Otherwise this role will pass to $(\tau_{a1},\tau_{a2})$. 

\smallskip

Finally, consider the case when $a=b+1$. In the previous
notation, if $\sigma_1^{(2)}\cap\tau_{a1}\ne \emptyset$,
then $(\tau_{a-1},\tau_{a1})$  is the bad pair we are looking
for, otherwise we should take $(\tau_{a1},\tau_{a2})$.

\bigskip

\centerline{{\it Subcase} (ii)(ii)${}^{\prime}$}

\medskip

Here $\sigma^{(1)}$ (resp. $\sigma^{(2)}$) breaks $\tau$
between $a$ and $a+1$ (resp. between $b$ and $b+1$),
and $a\ne b$.

\smallskip

If $a\ne b-1,\,b+1$, the modifications indiced in $\tau$ by
$\sigma^{(1)}$ and $\sigma^{(2)}$ do not interact
and the respective multiplications commute.

\smallskip

By symmetry, it remains to consider the case $a=b-1$. 
Choose $i\in\tau_a,\, j\in\tau_{a+1}.$ Summing first
over partitions $\alpha =(\tau_{a1},\tau_{a2})$
and $\beta =(\tau_{a+1,1},\tau_{a+1,2})$ we have
$$
l_{\sigma^{(1)}}\mu (\tau )=
-\sum_{\alpha :\,i\alpha}\mu (\dots \tau_{a1}\tau_{a2}\dots )
-\sum_{\beta :\,\beta j}\mu (\dots \tau_{a+1,1}\tau_{a+1,2}\dots ).
$$
Now, $\sigma^{(2)}$ will break the terms of the first (resp. second)
sum between $\tau_{a+1}$ and $\tau_{a+2}$ (resp. between
$\tau_{a+1,2}$  and $\tau_{a+2}$).
In order to multiply by $l_{\sigma^{(2)}}$ each term of these sums
we choose the same $j\in\tau_{a+1}$ and some $l\in\tau_{a+2}$.
Below we sum additionally over
2--partitions $\beta =(\tau_{a+1,1},\tau_{a+1,2})$
and $\gamma =(\tau_{a+2,1},\tau_{a+2,2})$
in the first two sums. In the second two sums the respective
notation is $\beta^{\prime}=(\tau_{a+1,2,1},\tau_{a+1,2,2})$:
$$
l_{\sigma^{(2)}}(l_{\sigma^{(1)}}\mu (\tau ))=
$$
$$
= \sum\Sb \alpha :\,i\alpha  \\ \beta :\,j\beta  \endSb
\mu (\dots \tau_{a1}\tau_{a2} \tau_{a+1,1}\tau_{a+1,2}\dots ) 
+\sum\Sb \alpha :\,i\alpha  \\ \gamma :\,\gamma l  \endSb
\mu (\dots \tau_{a1}\tau_{a2}\tau_{a+1} \tau_{a+2,1}\tau_{a+2,2}\dots )
$$
$$
+\sum\Sb \beta :\,\beta j  \\ \beta^{\prime} :\,j\beta^{\prime}  \endSb
\mu (\dots \tau_{a+1,1}\tau_{a+1,2,1} \tau_{a+1,2,2}\tau_{a+2}\dots ) 
+\sum\Sb \beta :\,\beta j  \\ \gamma :\,\gamma l  \endSb
\mu (\dots \tau_{a+1,1}\tau_{a+1,2} \tau_{a+2,1}\tau_{a+2,2}\dots )
\eqno(A.4)
$$
On the other hand, with the same notation we have:
$$
l_{\sigma^{(2)}}\mu (\tau )=
-\sum_{\beta :\,j\beta}\mu (\dots \tau_{a+1,1}\tau_{a+1,2}\dots )
-\sum_{\gamma :\,\gamma l}\mu (\dots \tau_{a+2,1}\tau_{a+2,2}\dots )
$$
and
$$
l_{\sigma^{(1)}}(l_{\sigma^{(2)}}\mu (\tau ))=
$$
$$
= \sum\Sb \alpha :\,i\alpha  \\ \beta :\,j\beta  \endSb
\mu (\dots \tau_{a1}\tau_{a2} \tau_{a+1,1}\tau_{a+1,2}\dots ) 
+\sum\Sb \beta :\, j\beta \\ \beta^{\prime\prime} :\,\beta^{\prime\prime} j  \endSb
\mu (\dots \tau_{a+1,1,1}\tau_{a+1,1,2} \tau_{a+1,2}\dots )
$$
$$
+\sum\Sb \alpha :\,i\alpha  \\ \gamma :\,\gamma l  \endSb
\mu (\dots \tau_{a1}\tau_{a2} \tau_{a+1}\tau_{a+2,1}\tau_{a+2,2}\dots ) 
+\sum\Sb \beta :\,\beta j  \\ \gamma :\,\gamma l  \endSb
\mu (\dots \tau_{a+1,1}\tau_{a+1,2} \tau_{a+2,1}\tau_{a+2,2}\dots )
\eqno(A.5)
$$
where $\beta^{\prime\prime}=(\tau_{a+1,1,1},\tau_{a+1,1,2})$.
Three of the four sums in (A.4) and (A.5) obviously
coincide. The third sum in (A.4) coincides with the
second sum in (A.5) because both are taken over
3--partitions of $\tau_{a+1}$ with $j$ in the middle part.

\bigskip

\centerline{{\it Subcase} (ii)(iii)${}^{\prime}$}

\medskip

Here $\sigma^{(1)}$  breaks $\tau$
between $a$ and $a+1$,  $\sigma^{(2)}$ does not break $\tau$.
We must check that $l_{\sigma^{(2)}}(l_{\sigma^{(1)}}\mu (\tau ))=0$,
by finding a bad pair for $\sigma^{(2)}$ and each term
in the right hand side of
$$
l_{\sigma^{(1)}}\mu (\tau )=
-\sum_{\alpha :\,i\alpha}\mu (\dots \tau_{a1}\tau_{a2}\dots )
-\sum_{\beta :\,\beta j}\mu (\dots \tau_{a+1,1}\tau_{a+1,2}\dots ).
$$
Denote by $(\tau_b,\tau_{b+1})$ a bad pair for
$\sigma^{(2)}$ and $\tau$. As in the subcase (i)(iii)${}^{\prime}$,
it will remain the bad pair unless $b\in \{a-1,a,a+1\},$
and will change somewhat in the exceptional cases.

\smallskip

More preciasely, if $b=a-1$, then for the terms of the second 
sum $(\tau_{a-1},\tau_{a})$ will be bad. For the first sum,
if $\sigma_1^{(2)}\cap \tau_{a1}\ne \emptyset$, the bad
pair will be $(\tau_{a-1},\tau_{a1})$. Otherwise it will
be $(\tau_{a1},\tau_{a2})$.

\smallskip

If $b=a$, then for the terms of the first 
sum $(\tau_{a2},\tau_{a+1})$ will be bad.  
For the second sum,
if $\sigma_1^{(2)}\cap \tau_{a+1,1}\ne \emptyset$, the bad
pair will be $(\tau_{a},\tau_{a+1,1})$. Otherwise it will
be $(\tau_{a+1,1},\tau_{a+1,2})$.

\smallskip

Finally, if $b=a+1$, then for the terms of the first 
sum $(\tau_{a+1},\tau_{a+2})$ will be bad.
For the second sum, it will be $(\tau_{a+1,2},\tau_{a+2})$.  

\medskip

In the last remaining case (iii)(iii)${}^{\prime}$ 
both multiplications produce zero.

\medskip

To complete the proof of the
Technical Lemma, it now remains to check that the elements (2.14), (2.15) generating
$I_B$ annihilate $H_{*B}$.

\medskip

{\bf A.4. Elements $r^{(1)}_{ij}$ annihilate $H_{*B}$.}

\smallskip

Fix $i,j$ and a partition $\tau$. If $\tau$ does not separate 
$i$ and $j$, we have $i,j\in\tau_a$ for some $a$, and then
$$
r^{(1)}_{ij}\mu (\tau ) =
\left(\sum_{\sigma :\,i\sigma j}l_{\sigma}-
\sum_{\sigma :\,j\sigma i}l_{\sigma}\right)\mu (\tau )
$$ 
$$
=\sum_{\alpha :\,i\alpha j}\mu (\tau (\alpha))-
\sum_{\alpha :\,j\alpha i}\mu (\tau (\alpha))
\eqno(A.6)
$$
where $\alpha$ runs over partitions of $\tau_a$. This expression vanishes
because of (2.23).

\smallskip

Assume now that $\tau$ separates $i$ and $j$, say, $i\in\tau_a,
j\in\tau_b,\,a<b.$ In this case $l_{\sigma}\mu (\tau )=0$
for all $\sigma$ with $j\sigma i.$ The remaining terms of (A.6)
vanish unless $\sigma$ breaks $\tau$ at some $\tau_c,\,a\le c\le b$,
or else between $\tau_c$ and $\tau_{c+1}$ for $a\le c \le b-1.$
In the latter cases each term corresponding to one $\sigma$
can be replaced by a sum of terms corresponding to
the 2--partitions $\alpha_c$ of $\tau_c$ with the help of (2.24) and (2.25).

\smallskip

Let us choose $k_c\in\tau_c$ for all $a\le c\le b$ so that
$k_a=i,\,k_b=j$ and spell out the resulting expression:
$$
\sum_{\sigma :\,i\sigma j}l_{\sigma}\mu (\tau )=
\sum_{c=a}^{c=b}{}^{\prime}\left(\sum_{\alpha_c:\,k_c\alpha_c}+
\sum_{\alpha_c:\,\alpha_c k_c}\right) \mu (\tau (\alpha_c ))
$$
$$
-\sum_{c=a}^{c=b-1}\left(\sum_{\alpha_c:\,k_c\alpha_c}
+\sum_{\alpha_{c+1}:\,\alpha_{c+1}k_{c+1}}\right) \mu (\tau (\alpha_{c+1} ))\, .
$$
Here prime at the first sum indicates that the terms with $\alpha_ai$
and $j\alpha_b$ should be skipped.

\smallskip

All the terms of this expression cancel.

\medskip

{\bf A.5. Elements $r^{(2)}(\sigma^{(1)},\sigma^{(2)})$ annihilate $H_{*B}$.}

\smallskip

These elements correspond
to the pairs  ($\sigma^{(1)},\,\sigma^{(2)}$) that do not break each other.
If at least one of them, say $\sigma^{(1)}$, does not break $\tau$ either, then
$l_{\sigma^{(1)}}\mu (\tau )=0$ so that 
$r^{(2)}(\sigma^{(1)},\sigma^{(2)})\mu (\tau )=0.$ 
If both $\sigma^{(1)},\,\sigma^{(2)}$ break $\tau$, a contemplation
will convince the reader that they must break $\tau$ at one and the same
component $\tau_a.$ This is the subcase (i)(i)${}^{\prime}$ of A.3,
and we made the relevant comment there.

\newpage

\centerline{\bf References}

\medskip

[Da] V.~Danilov. {\it The geometry of toric varieties.} Russian
Math. Surveys, 33:2 (1978), 97--54.

\smallskip

[Du1] B.~Dubrovin. {\it Geometry of 2D topological fielld theories.}
In: Springer LNM, 1620 (1996), 120--348.

\smallskip

[Du2] B.~Dubrovin. {\it Painlev\'e transcendents in two--dimensional
topological field theory.} Preprint math.AG/9803107.

\smallskip

[Du3] B.~Dubrovin. {\it Flat pencils of metrics and Frobenius manifolds.}
Preprint math.AG/9803106

\smallskip

[Fu] W.~Fulton. {\it Introduction to toric varieties.}
Ann. Math. Studies, Nr. 131, Princeton University Press,
Princeton NJ, 1993.

\smallskip

[GeSe] I.~M.~Gelfand, V.~Serganova. {\it Combinatorial geometries
and torus strata on compact homogeneous spaces.} Uspekhi Mat. Nauk
42 (1987), 107--134.

\smallskip

[Ka1]  M.~Kapranov. {\it The permutoassociahedron, MacLane's 
coherence theorem and asymptotic zones for the KZ equation.}
Journ. of Pure and Appl. Algebra, 83 (1993), 119--142.

\smallskip

[Ka2] M.~Kapranov. {\it Chow quotients of Grassmannians. I.}
Advances in Soviet Math., 16:2 (1993), 29--110.

\smallskip

[Ke] S.~Keel. {\it Intersection theory of moduli spaces
of stable $n$--pointed curves of genus zero.} Trans. AMS,
330 (1992), 545--574.

\smallskip

[Kn1] F.~Knudsen. {\it Projectivity of the moduli space of stable
curves, II: the stacks $M_{g,n}.$} Math. Scand. 52 (1983),
161--199.

\smallskip

[Kn2] F.~Knudsen. {\it The projectivity of the moduli space
of stable curves III: The line bundles on $M_{g,n}$ and a proof
of projectivity of $\overline{M}_{g,n}$ in characteristic 0.}
Math. Scand. 52 (1983), 200--212. 

\smallskip

[KM1] M.~Kontsevich, Yu.~Manin. {\it Gromov--Witten classes, quantum cohomology, and enumerative geometry.} Comm. Math. Phys.,
164:3 (1994), 525--562.

\smallskip

[KM2]  M.~Kontsevich, Yu.~Manin. {\it Relations between the correlators of the
topological sigma--model coupled to gravity.} Comm. Math. Phys.,
196 (1998), 385--398.

\smallskip

[KMK] M.~Kontsevich, Yu.~Manin. {\it Quantum cohomology of a product (with Appendix by R. Kaufmann)}. Inv. Math., 124, f. 1--3 (1996), 313--339.

\smallskip 

[Lo1] A.~Losev. {\it Commutativity equations, operator--valued
cohomology of the ``sausage'' compactification of
$(\bold{C}^*)^N/\bold{C}^*$ and SQM.} Preprint ITEP--TH--84/98,
LPTHE--61/98.

\smallskip

[Lo2] A.~Losev. {\it Passing from solutions to Commutativity
Equations to solutions to Associativity Equations and background
independence for gravitational descendents.} Preprint ITEP--TH--85/98,
LPTHE--62/98.

\smallskip

[Ma] Yu.~I.~Manin. {\it Frobenius manifolds, quantum cohomology, and moduli
spaces.} AMS Colloquium Publications, vol. 47, Providence, RI, 1999,
xiii+303 pp.

\smallskip

[MaZ] Yu.~I.~Manin, P.~Zograf. {\it Invertible Cohomological Field Theories and Weil--Petersson volumes.} 
Preprint math.AG/9902051 (to appear in Ann. Inst. Fourier).

\smallskip

[Si] C.~Simpson. {\it The Hodge filtration on non--abelian
cohomology.} Preprint, 1996.

\enddocument